\documentclass[pdflatex,sn-mathphys-num]{sn-jnl}

\usepackage{amsmath,amssymb,mathrsfs}
\usepackage{listings}
\usepackage{xcolor}
\usepackage{tikz-cd}

\theoremstyle{plain}
\newtheorem{theorem}{Theorem}
\newtheorem{proposition}[theorem]{Proposition}
\newtheorem{lemma}[theorem]{Lemma}
\newtheorem{corollary}[theorem]{Corollary}

\theoremstyle{definition}
\newtheorem{definition}[theorem]{Definition}
\newtheorem{remark}[theorem]{Remark}

\newcommand{\Spec}{\operatorname{Spec}}
\newcommand{\Proj}{\operatorname{Proj}}
\newcommand{\Rees}{\operatorname{Rees}}
\newcommand{\Pot}{\operatorname{Pot}}
\newcommand{\Clos}{\operatorname{Clos}}

\newcommand{\Bl}{\operatorname{Bl}}
\newcommand{\BlM}{\operatorname{Bl}^{\mathrm{aff}}}
\newcommand{\Pres}{\mathcal{M}}

\newcommand{\Nat}{\mathbb{N}}
\newcommand{\Int}{\mathbb{Z}}
\newcommand{\cL}{\mathcal{L}}
\newcommand{\cI}{\mathcal{I}}
\newcommand{\cZ}{\mathcal{Z}}
\newcommand{\cU}{\mathcal{U}}
\newcommand{\cC}{\mathcal{C}}
\newcommand{\sO}{\mathcal{O}}

\definecolor{leankw}{rgb}{0.15,0.25,0.65}
\definecolor{leanbg}{rgb}{0.975,0.975,0.98}
\lstdefinestyle{lean}{
  basicstyle=\ttfamily\footnotesize,
  keywordstyle=\color{leankw}\bfseries,
  backgroundcolor=\color{leanbg},
  frame=single,
  rulecolor=\color{gray!40},
  columns=fullflexible,
  breaklines=true,
  breakatwhitespace=true,
  showstringspaces=false,
  xleftmargin=8pt,
  xrightmargin=4pt,
  aboveskip=8pt,
  belowskip=8pt,
  morekeywords={def,lemma,theorem,structure,instance,abbrev,noncomputable,
    where,with,fun,by,exact,intro,refine,have,let,rw,simp,apply,variable,
    open,namespace,end,extends,Prop,Type,Set,Sort,do,then,else,if,match},
  literate=
    {α}{{$\alpha$}}1 {β}{{$\beta$}}1 {γ}{{$\gamma$}}1 {δ}{{$\delta$}}1
    {ε}{{$\varepsilon$}}1 {ι}{{$\iota$}}1 {κ}{{$\kappa$}}1 {λ}{{$\lambda$}}1
    {μ}{{$\mu$}}1 {ν}{{$\nu$}}1 {ρ}{{$\rho$}}1 {σ}{{$\sigma$}}1
    {τ}{{$\tau$}}1 {φ}{{$\varphi$}}1 {χ}{{$\chi$}}1 {ψ}{{$\psi$}}1
    {ω}{{$\omega$}}1 {Ψ}{{$\Psi$}}1 {Σ}{{$\Sigma$}}1 {Π}{{$\Pi$}}1
    {ℕ}{{$\mathbb{N}$}}1 {ℤ}{{$\mathbb{Z}$}}1 {ℚ}{{$\mathbb{Q}$}}1
    {𝒜}{{$\mathcal{A}$}}1 {ℱ}{{$\mathcal{F}$}}1 {𝒰}{{$\mathcal{U}$}}1
    {𝒱}{{$\mathcal{V}$}}1 {𝒬}{{$\mathcal{Q}$}}1 {𝟙}{{$\mathbb{1}$}}1
    {→}{{$\rightarrow$}}1 {↦}{{$\mapsto$}}1 {←}{{$\leftarrow$}}1
    {⟶}{{$\longrightarrow$}}2 {⟹}{{$\Longrightarrow$}}2
    {≅}{{$\cong$}}1 {≃}{{$\simeq$}}1 {≠}{{$\neq$}}1 {≤}{{$\leq$}}1
    {≥}{{$\geq$}}1 {≈}{{$\approx$}}1 {≡}{{$\equiv$}}1
    {∈}{{$\in$}}1 {∉}{{$\notin$}}1 {∀}{{$\forall$}}1 {∃}{{$\exists$}}1
    {∧}{{$\wedge$}}1 {∨}{{$\vee$}}1 {¬}{{$\neg$}}1 {∅}{{$\emptyset$}}1
    {⟨}{{$\langle$}}1 {⟩}{{$\rangle$}}1 {≪}{{$\ll$}}1 {≫}{{$\gg$}}1
    {⧸}{{$/$}}1 {↘}{{$\searrow$}}1 {∘}{{$\circ$}}1 {•}{{$\cdot$}}1
    {∩}{{$\cap$}}1 {∪}{{$\cup$}}1 {⊆}{{$\subseteq$}}1 {⊂}{{$\subset$}}1
    {⊤}{{$\top$}}1 {⊥}{{$\bot$}}1 {∑}{{$\sum$}}1 {⨁}{{$\bigoplus$}}1
    {∣}{{$\mid$}}1 {⇑}{{}}0 {₀}{{$_0$}}1 {₁}{{$_1$}}1 {₂}{{$_2$}}1
    {ₐ}{{$_a$}}1 {ᵥ}{{$_v$}}1 {⁻¹}{{$^{-1}$}}2 {·}{{$\cdot$}}1
    {√}{{$\sqrt{}$}}1 {∼}{{$\sim$}}1 {Γ}{{$\Gamma$}}1
}
\newcommand{\leansrc}[2]{%
  \par\vspace{-4pt}\noindent{\scriptsize\ttfamily #1, lines #2}\par\vspace{2pt}}

\begin{document}

\title{Projective blowups: a formal and multicentered proof}
\author[1]{Arnaud Mayeux}
\affil[1]{University of Wisconsin-Madison, Madison, WI, USA. \newline\email{mayeux@wisc.edu}}

\abstract{
We give a machine-checked proof, in Lean 4, of the universal property of
multicentered blowups: given a finite family $Z$ of closed subschemes of a scheme
$X$, presented on each chart by ideals, there is a scheme $\Bl_Z X$ over
$X$, terminal among $X$-schemes on which every member of the family becomes
an effective Cartier divisor. To the best of our knowledge this is the first
formalisation of blowups, single- or multicentered, in any theorem prover.
The proof is organised around the universal property of the multicentered
dilatation of a ring. The comparison isomorphisms on
overlaps, the cocycle identity, and the independence of the chosen
presentation are obtained from the uniqueness part of this universal
property. We give the exact theorem
statement as formalized, the proof architecture, and the main new
mathematical ingredients, together with the corresponding Lean code.
}


\maketitle

\tableofcontents

\section{Introduction}\label{sec:intro}

Blowing up a scheme along a closed subscheme is one of the foundational operations
in algebraic geometry. The classical blowup solves a lifting problem: any $X$-scheme
on which the centre becomes a Cartier divisor admits a unique morphism to the blowup.

For a finite family of centres, the multicentered blowup is the universal space on which
all the centres become Cartier simultaneously. This paper establishes the universal property
of multicentered blowups with a machine-checked proof in Lean 4 \cite{lean4}. To the best of our
knowledge this is the first formalisation of blowups, single- or multicentered, in any
theorem prover.

The key tool is multi-centered dilatations and their universal property, which turns ideals into principal ideals. The global blowup is
obtained by gluing dilatation charts using the multigraded $\Proj$ construction of
Brenner and Schröer.

The supporting machinery had to be developed. The formalization of the multigraded $\Proj$ and multicentered dilatations were established in prior work \cite{mayeux-zhang,mayeux-zhang-short}. The result is a complete, checked proof of
a non-trivial theorem in algebraic geometry, with no noetherian, properness, or finite
type assumptions on $X$, and all dependences on the ideals explicit.

This document is organised as follows. Section 1.1 situates the work relative to prior
results and existing libraries. Section 1.2 outlines the proof strategy, emphasising the
role of the universal property of dilatations. Sections 2--4 contain the proof: the theorem
statement, the dilatation and its properties, the reduction to the affine case, and
the globalisation. Section 5 lists the correspondence between the statements
of this paper and their Lean declarations.

\subsection{Relation to other work}\label{sec:relation}

The theorem proved below is not new as a statement, and neither is the
dilatation on which its proof rests. The dilatation comes from
\cite{mayeux-dilatations}, the published version, where multicentered
dilatations of rings, schemes and algebraic spaces are introduced, and where
the properties used here as Theorem \ref{thm:dil-univ} and Lemmas
\ref{lem:dil-nzd} and \ref{lem:dil-gen} are established. The multicentered
blowup and its universal property, together with the use of the multigraded
$\Proj$ of Brenner and Schr\"oer \cite{brenner-schroer} to glue the dilatation
charts into a scheme, appeared in the fourth arXiv version of that paper
\cite{mayeux-dilatations-v4}, but were dropped from the published version of
\cite{mayeux-dilatations}. A second proof, next to a detailed Brenner–Schröer
theory, was given in \cite{mayeux-riche}. For $X = \Spec A$ with closed subschemes $Z_1, \dots, Z_k$ cut out by ideals
$L_1, \dots, L_k$, the affine blowup is defined as
\[
  \Bl_Z X \;:=\; \Proj\bigl( \Bl_{L_1, \dots, L_k} A \bigr),
\]
the multi-graded $\Proj$ of the multi-Rees algebra. It is a final object of the category
$\mathrm{Sch}_X^{Z_1, \dots, Z_k\text{-reg}}$ of $X$-schemes on which all the
$Z_i$ become regularly immersed. That is the affine case of Corollary
\ref{cor:terminal}, proved below as Theorem \ref{thm:aff-main}.

Both ingredients matter, and they enter on an equal footing: the dilatations
are what the charts are, and the multi-graded $\Proj$ is what makes them into a
scheme. For the latter, \cite{mayeux-dilatations-v4} refers to the
multihomogeneous spectra of Brenner and Schr\"oer \cite{brenner-schroer}, whose
treatment is systematically developed in \cite{mayeux-riche}. In \cite{mayeux-riche} the statement is proved for
general $X$, as \cite[Proposition 3.30]{mayeux-riche}: if $T \to X$ is such
that every $T \times_X Y_i$ is an effective Cartier divisor on $T$, then there
is a unique morphism of $X$-schemes $T \to X'$.

The formalisation of the Brenner–Schröer multi-graded $\Proj$ construction that
the present work relies on throughout was developed by Mayeux and Zhang in
\cite{mayeux-zhang} (see also the short account \cite{mayeux-zhang-short}),
following the algebraic treatment of \cite{mayeux-riche}. That paper formalizes multi-graded $\Proj$
schemes and dilatations of rings in Lean 4, as local, ring-level
constructions; it stops short of gluing dilatations of quasi-coherent
algebras into a global blowup of schemes, a limitation the authors note
explicitly. The present work builds on that formalisation without mathematical
modification, using it as a black box for everything concerning potions,
relevant homogeneous submonoids, and the gluing of $\Proj$ from good potion
ingredients, and develops the global multicentered blowup and establishes its universal
property. Mathlib itself contains no theory of blowups, single- or
multicentered, and the multicentered blowup is constructed entirely from these primitives and scheme-theoretic foundations in mathlib.

The construction synthesizes the multi-graded $\Proj$ and dilatations as complementary
components: dilatations define the affine charts, and the multi-graded $\Proj$
provides the gluing machinery. The affine case identifies the charts
$D_+(f)$ of $\Bl_Z X$ with dilatations $A[L_1/f, \dots, L_k/f]$, then
lifts the local universal property of dilatations to the global setting via the
$\Proj$ construction. Sections \ref{sec:dilatation} and
\ref{ssec:cars} to \ref{ssec:aff-exist} develop this argument in the
form the formalisation uses.

The present account provides the globalisation of the universal property and its complete proof. The
proof of Theorem \ref{thm:main} is organised so that every step beyond the
affine case, namely the base change along open immersions, the comparison
isomorphisms on overlaps, the triple overlap maps, the cocycle identity, and
the independence of the chosen presentation, is deduced from the affine
universal property by uniqueness. The bookkeeping of covers and defining ideals that this requires is
managed by the notions of presented centre and centre of Section
\ref{ssec:pullback}.
Multi-graded $\Proj$ and blowups of this kind are used elsewhere. Tang
\cite{tang} constructs Gysin maps for regular immersions in a
$\mathbf{P}^1$-unstable, non-$\mathbf{A}^1$-invariant theory of motivic spaces
and spectra, using pushout-blowups. About further functoriality of these blowups, he writes
\cite[Remark~A.32]{tang}: ``It will probably involve the $\Proj$ construction of
$\Nat^n$-graded rings, i.e.\ an animated version of the theory in''
\cite{brenner-schroer} and \cite{mayeux-riche}.

The formalisation is built on mathlib \cite{mathlib}, which supplies schemes,
$\Spec$, fibre products, open and closed immersions, gluing data, homogeneous
localisation and the theory of non-zero-divisors. The general scheme theory in mathlib, on which everything here depends, was developed in \cite{buzzard-schemes}. Building on this, Zhang \cite{zhang-proj} formalized the monograded Proj construction, providing the first non-affine scheme formalization in a theorem prover.

Beyond this general infrastructure, mathlib does not yet contain blowups,
dilatations, multigraded $\Proj$, or effective Cartier ideal sheaves. It
defines closed immersions, quasi-coherent ideal sheaves
(\texttt{IdealSheafData}), and the correspondence between them
(\texttt{overEquivIdealSheafData}). We present a centre by ideals of the rings
of an affine cover; Subsection \ref{ssec:mathlib} shows that this presentation
agrees with mathlib's ideal sheaves. We recall the classical statement.

\begin{lemma}[{\cite[\href{https://stacks.math.columbia.edu/tag/01QN}{Tag
01QN}]{stacks-project}}]\label{lem:closed-imm}
Let $i : Z \to X$ be a morphism of schemes. The following are equivalent.
\begin{enumerate}
\item $i$ is a closed immersion.
\item For every affine open $U = \Spec R \subseteq X$ there is an ideal
      $I \subseteq R$ with $i^{-1}(U) \cong \Spec(R/I)$ as schemes over $U$.
\item There is an affine open cover $X = \bigcup_{j \in J} U_j$, with
      $U_j = \Spec R_j$, and for each $j$ an ideal $I_j \subseteq R_j$ with
      $i^{-1}(U_j) \cong \Spec(R_j / I_j)$ as schemes over $U_j$.
\item $i$ induces a homeomorphism of $Z$ onto a closed subset of $X$ and
      $i^{\sharp} : \sO_X \to i_* \sO_Z$ is surjective.
\item $i$ induces a homeomorphism of $Z$ onto a closed subset of $X$,
      $i^{\sharp}$ is surjective, and $\ker(i^{\sharp}) \subseteq \sO_X$ is a
      quasi-coherent sheaf of ideals.
\item $i$ induces a homeomorphism of $Z$ onto a closed subset of $X$,
      $i^{\sharp}$ is surjective, and $\ker(i^{\sharp}) \subseteq \sO_X$ is a
      sheaf of ideals which is locally generated by sections.
\end{enumerate}
\end{lemma}

Condition (3) is the one we use. Definition \ref{def:preclos} below is that
condition with the cover and the ideals promoted from existentially quantified
objects to data: rather than asserting that suitable $U_j$ and $I_j$ exist, a
presented centre carries a specific choice of them. This is what makes the
Cartier condition of Definition \ref{def:cartier} expressible, since that
condition is a statement about the chosen ideals, and it is why the equivalence
relation of Definition \ref{def:rel} has to be imposed afterwards to forget the
choice again.

In the formalisation, a presented centre supplies the data of (3), with charts
$\Spec R_j \to X$. It is proved that each $Z_i \to X$ is then a closed
immersion, and that its kernel restricts on each chart to $I_j$
(Subsection \ref{ssec:mathlib}).

\subsection{Strategy of the proof}\label{sec:strategy}

Roughly, the proof covers $X$ by affine charts, treats each chart, and glues.
On a chart, the main input is the universal property of the multicentered
dilatation; the gluing relies mostly on the uniqueness part of the affine
universal property.

The proof does not follow \cite{mayeux-riche} closely, nor \cite{mayeux-dilatations-v4}. The argument there uses
multicentered global dilatations of schemes and rests on a body of supporting
theory: multi-graded $\Proj$ and its charts, quasi-coherent sheaves, closed
subschemes, Cartier divisors and their behaviour under base change. Most of that
is not available in mathlib, as recorded in Section \ref{sec:relation}. The
present proof operates at the ring level, avoiding the need for global scheme-theoretic infrastructure, and constructs the global blowup
via covers, charts, gluing data, and uniqueness arguments. Following
\cite{mayeux-riche} literally would have meant developing both the global
dilatations of schemes and all of their supporting theory, which is less
elementary than the road taken here via dilatations of rings. The multigraded
$\Proj$ construction of Brenner and Schröer \cite{brenner-schroer} is used to
define the blowup by gluing dilatation charts. The universal property of the
multicentered dilatation of a ring is then used to prove the universal property
of the projective blowup locally, chart by chart. Every global step
beyond that, including the gluing and comparison isomorphisms, uses the
uniqueness part of the local universal property of the multicentered dilatation.

The shape of the argument is the following.

\emph{Existence, locally.} In the affine local situation the existence half of
the universal property of dilatations applies. Its hypothesis is the Cartier
condition with prescribed generators: the images of the $a_j$ are
non-zero-divisors and generate the images of the large ideals. So
the required morphism into a chart is produced by Theorem \ref{thm:dil-univ}.

\emph{Gluing, by uniqueness.} The local morphisms are then glued using the
uniqueness part of the same universal property: two local solutions solve the
same problem on an overlap, hence agree there. The same
argument gives the remaining steps: the comparison
isomorphisms $t_{\gamma\delta}$ on double overlaps, the triple overlap maps
$t'_{\gamma\delta\epsilon}$, the cocycle identity, the base change isomorphism
of Proposition \ref{prop:bc}, and the independence of the chosen presentation.
In each case two objects are shown to solve one problem, and are therefore
canonically identified.

Sections \ref{sec:dilatation} to \ref{sec:proof} carry this out.

\paragraph{Formal artifact.} The Lean code accompanying this paper is
publicly available at \url{https://github.com/rndmx/MULTI-PROJ-BLOW} \cite{blowups-code}, version
\texttt{v2.2}. It builds with Lean \texttt{v4.22.0-rc3} and mathlib at commit
\texttt{7b332d3}; running \texttt{lake exe cache get} followed by
\texttt{lake build} checks all proofs. The formalization contains no \texttt{sorry}; the main
theorems depend only on the axioms \texttt{propext}, \texttt{Classical.choice}
and \texttt{Quot.sound}. Every code excerpt in this paper is
preceded by its file and line numbers in that version.
The repository includes portions of the code of \cite{mayeux-zhang}, updated to
this version of mathlib, so that it builds on its own.

\section{The theorem}

\subsection{Centres}

Throughout, $X$ is a scheme. The object we blow up is a family of closed
subschemes, together with an affine cover on whose charts each member is cut
out by an ideal. The definitions allow an arbitrary index type; finiteness is
assumed in the theorems. For infinite families, a blowup need not exist: there exist a scheme $X$ and an infinite family of closed subschemes of $X$ that has no blowup.

\begin{definition}\label{def:preclos}
A \emph{presented centre} on $X$ consists of:
\begin{enumerate}
\item a type $\iota$ of indices;
\item for each $i \in \iota$, a scheme $Z_i$ together with a morphism
      $Z_i \to X$;
\item an affine open cover $\cU = \bigl( u_\gamma : U_\gamma := \Spec A_\gamma \to X
      \bigr)_{\gamma \in J}$ of $X$ by open immersions;
\item for each $i$ and $\gamma$, an ideal $I_{i\gamma} \subseteq A_\gamma$;
\item for each $i$ and $\gamma$, an isomorphism of $\Spec A_\gamma$-schemes
\[
  \Spec\bigl(A_\gamma / I_{i\gamma}\bigr)
  \;\xrightarrow{\ \sim\ }\;
  Z_i \times_X \Spec A_\gamma .
\]
\end{enumerate}
\end{definition}

Condition (5) does two things at once. It forces $Z_i \times_X U_\gamma$ to be
affine, and it says that $I_{i\gamma}$ is the ideal cutting it out. So $Z_i$ is
a closed subscheme of $X$, presented on each chart. In Lean this is a single
structure.

\leansrc{Project/Blowups/PreClosAndClos.lean}{52--64}
\begin{lstlisting}
variable (X) in
structure PreClos where
  (indnumb : Type)
  (subscheme : indnumb → Scheme)
  [over : ∀ (i : indnumb), Scheme.Over (subscheme i) X]
  cov : Scheme.AffineCover.{u+1} (P := @IsOpenImmersion) X
  ideal : ∀ (_ : indnumb) (γ : cov.J), Ideal (cov.obj γ)
  condiso : ∀ (i : indnumb) (γ : cov.J),
    Spec (CommRingCat.of (cov.obj γ ⧸ ideal i γ)) ≅
    pullback (f := subscheme i ↘ X) (g := cov.map γ)
  condover : ∀ (i : indnumb) (γ : cov.J),
    Scheme.Hom.IsOver (condiso i γ).hom
      (Spec (CommRingCat.of (cov.obj γ)))
\end{lstlisting}

The cover is auxiliary. Two presented centres should be regarded as the same
centre when they cut out the same subschemes, whatever charts they use to do
it. This is the relation we quotient by.

\begin{definition}\label{def:rel}
Two presented centres $Z, Z'$ are \emph{equivalent} if there is a bijection
$\sigma : \iota \to \iota'$ of index types and isomorphisms
$Z_i \cong Z'_{\sigma(i)}$ over $X$ for every $i$. A \emph{centre} on $X$ is an
equivalence class. We write $\Clos(X)$ for the set of centres.
\end{definition}

\leansrc{Project/Blowups/PreClosAndClos.lean}{73--76, 111, 123--124}
\begin{lstlisting}
structure relStructure (Z Z' : PreClos X) where
  indnumb_equiv : Z.indnumb ≃ Z'.indnumb
  subscheme_iso : ∀ i, Z.subscheme i ≅ Z'.subscheme (indnumb_equiv i)
  subscheme_iso_over : ∀ i, Scheme.Hom.IsOver (subscheme_iso i).hom X

def rel : PreClos X → PreClos X → Prop := fun Z Z' => Nonempty (relStructure Z Z')

variable (X)
def Clos := Quotient (relSetoid X)
\end{lstlisting}

The covers play no role in the relation. That is the whole point: a centre is a
family of closed subschemes of $X$, and the chart data used to present it is
forgotten. Everything we build will be built from a presentation, and a
substantial part of the work consists in showing that the result does not
depend on the presentation.

\subsection{Pullback of centres}\label{ssec:pullback}

Centres pull back along an arbitrary morphism, and the construction is the
obvious one, but it is worth writing out because the choice of cover on the
target is forced and reappears later.

Let $f : X' \to X$ be a morphism and $Z$ a presented centre on $X$ with cover
$(u_\gamma : \Spec A_\gamma \to X)_{\gamma \in J}$. The pulled back centre
$f^*Z$ has:
\begin{itemize}
\item the same index type $\iota$, with subschemes $Z_i \times_X X'$;
\item as cover, the family indexed by pairs $(\gamma, \beta)$ where $\gamma \in
      J$ and $\beta$ runs over an affine open cover of $X' \times_X U_\gamma$,
      with chart map
      \[
        \Spec A_{\gamma\beta} \longrightarrow X' \times_X U_\gamma
        \xrightarrow{\ \mathrm{pr}_1\ } X';
      \]
\item as ideal at $(i, (\gamma,\beta))$, the image
      $I_{i\gamma} \cdot A_{\gamma\beta}$ under the induced ring map
      $A_\gamma \to A_{\gamma\beta}$.
\end{itemize}
The point is that one cannot simply pull back the cover: $X' \times_X U_\gamma$
is open in $X'$ but need not be affine, so it has to be covered again. This
double indexing is why the index type of a pulled back centre is a sigma type,
and it is the reason several later statements quantify over pairs.

That $f^*$ respects the equivalence relation of Definition \ref{def:rel} is
straightforward, since the relation only involves the subschemes $Z_i$ and
base change of schemes preserves isomorphism over the base. So $f^*$ descends
to a map
\[
  f^* : \Clos(X) \longrightarrow \Clos(X').
\]
Functoriality, $(g \circ f)^* = f^* \circ g^*$, holds at the level of $\Clos$
though not on the nose at the level of presentations, because the chosen affine
covers differ. This is a first instance of a phenomenon that recurs throughout:
statements that are false for presentations become true for centres.

\leansrc{Project/Blowups/PreClosAndClos.lean}{351--357, 389--390}
\begin{lstlisting}
def pullback_PreClos (X' : Scheme) (f : X' ⟶ X) (Z : PreClos X) : PreClos X' where
  indnumb := Z.indnumb
  subscheme i := pullback (Z.subscheme i ↘ X) f
  over i := ⟨pullback.snd _ _⟩
  cov := pull_cov X Z X' f
  ideal i γβ := pull_ideal X Z X' f γβ i
  condiso i γβ := pullback_PreClos_condiso _ γβ i

def pullback_Clos {X' : Scheme} (f : X' ⟶ X) : Clos X → Clos X' :=
  Quotient.map (pullback_PreClos X X' f) <| fun Z Z' e => Nonempty.map (pullback_lem X Z Z' X' f) e
\end{lstlisting}

\subsection{Cartier centres}

The following condition is stated on the charts of a presentation. For a centre
given by a finite family of ideal sheaves, it is equivalent to each ideal sheaf
being effective Cartier (Proposition \ref{prop:cartier-mathlib}).

\begin{definition}\label{def:cartier}
A presented centre $Z$ is \emph{Cartier} if for every $i \in \iota$ and every
$\gamma \in J$ the ideal $I_{i\gamma} \subseteq A_\gamma$ is principal and
generated by a non-zero-divisor. A centre $\cZ \in \Clos(X)$ is Cartier if
\emph{some} presentation of $\cZ$ is Cartier.
\end{definition}

\leansrc{Project/Blowups/PreClosAndClos.lean}{126--146}
\begin{lstlisting}
structure PrePri extends PreClos X where
  [prin : ∀ i γ, ideal i γ |>.IsPrincipal]

structure IsPrePri (Z : PreClos X) : Prop where
  prin : ∀ i γ, Z.ideal i γ |>.IsPrincipal

attribute [instance] PrePri.prin

structure PreCars extends PrePri X where
  nonzerodiv : ∀ i γ, Submodule.IsPrincipal.generator (ideal i γ) ∈ nonZeroDivisors (cov.obj γ)

structure IsPreCars (Z : PreClos X) extends IsPrePri _ Z where
  nonzerodiv : ∀ i γ, Submodule.IsPrincipal.generator (Z.ideal i γ) ∈ nonZeroDivisors (Z.cov.obj γ)

structure IsPri (Z : Clos X) : Prop where
  exists_rep : ∃ (Z' : PreClos X), IsPrePri _ Z' ∧ Quotient.mk'' Z' = Z

def Pri : Set (Clos X) := {x : Clos X | IsPri _ x}

structure IsCars (Z : Clos X) : Prop where
  exists_rep : ∃ (Z' : PreClos X), IsPreCars _ Z' ∧ Quotient.mk'' Z' = Z
\end{lstlisting}

Two features of Definition \ref{def:cartier} deserve emphasis, because they are
what make the proof work and they are easy to blur in informal writing.

First, the existential quantifier sits outside. A centre is Cartier when there
is \emph{one} presentation on which all the ideals are principal on all charts.
It is not required that every presentation have this property, and indeed most
do not.

Second, the generator is required to be a non-zero-divisor. Together with
principality this says exactly that each $Z_i$ pulls back to an effective
Cartier divisor.

\subsection{Statement}

Fix a presented centre $Z$ on $X$ with $\iota$ finite. We will construct a
scheme $\Bl_Z X$ over $X$. The theorem characterising it is the following.

\begin{theorem}[Universal property, presented form]\label{thm:main-pre}
Let $Z$ be a presented centre on $X$ with finite index type. Let $T$ be a
scheme over $X$, with structure morphism $f : T \to X$. If the pulled back
presented centre $f^*Z$ is Cartier, then there is exactly one morphism
\[
  \varphi : T \longrightarrow \Bl_Z X
\]
over $X$.
\end{theorem}

In Lean, with $\Bl_Z X$ written \texttt{BlGlob Z}:

\leansrc{Project/Blowups/BlowupExists.lean}{3009--3016}
\begin{lstlisting}
lemma PreProjBlowup_UnivProp
    (Z : PreClos X)
    [DecidableEq Z.indnumb]
    [Fintype Z.indnumb]
    [(i : Z.indnumb →₀ ℤ) → Decidable (i ∈ Set.range (ρNatToInt Z.indnumb))]
    {T : Scheme} [T.Over X]
    (cond : IsPreCars _ <| pullback_PreClos _ _ (T ↘ X) Z) :
    ∃! φ : T ⟶ BlGlob Z, Scheme.Hom.IsOver φ X := by
\end{lstlisting}

For $\cZ \in \Clos(X)$ with presentation $Z$, set $\Bl_\cZ X := \Bl_Z X$
(well defined by Proposition \ref{prop:rel}).

\begin{theorem}[Universal property]\label{thm:main}
Let $\cZ \in \Clos(X)$ be a centre with finite index type and let $T$ be a
scheme over $X$ with structure morphism $f$. If $f^*\cZ$ is Cartier, then there
is exactly one morphism $T \to \Bl_\cZ X$ over $X$.
\end{theorem}

\leansrc{Project/Blowups/BlowupExists.lean}{3395--3398}
\begin{lstlisting}
theorem GlobalBlowup_UnivProp (Z : Clos X) [Finite Z.out.indnumb]
    {T : Scheme} [T.Over X]
    (cond : IsCars _ <| pullback_Clos (T ↘ X) Z) :
    ∃! φ : T ⟶ GlobalBlowup Z, Scheme.Hom.IsOver φ X := by
\end{lstlisting}

Alongside these sits the statement that the blowup itself has the property in
question.

\begin{theorem}\label{thm:blglob-cartier}
The centre pulled back to $\Bl_Z X$ is Cartier.
\end{theorem}

\leansrc{Project/Blowups/BlowupExists.lean}{3342--3345}
\begin{lstlisting}
lemma BlGlob_IsCars (Z : PreClos X)
    [DecidableEq Z.indnumb] [Fintype Z.indnumb]
    [(i : Z.indnumb →₀ ℤ) → Decidable (i ∈ Set.range (ρNatToInt Z.indnumb))] :
    IsCars _ (pullback_Clos ((BlGlob Z) ↘ X) (Quotient.mk'' Z)) := by
\end{lstlisting}

The scheme $\Bl_\cZ X$ itself, for a centre rather than a presentation, is
obtained by choosing the representative that the quotient supplies.

\leansrc{Project/Blowups/BlowupExists.lean}{3382--3386}
\begin{lstlisting}
noncomputable def GlobalBlowup (Z : Clos X) [Finite Z.out.indnumb] : Scheme :=
  letI := Fintype.ofFinite Z.out.indnumb
  letI := Classical.decEq Z.out.indnumb
  letI := Classical.decPred (· ∈ Set.range (ρNatToInt Z.out.indnumb))
  BlGlob Z.out
\end{lstlisting}

Combining Theorems \ref{thm:main} and \ref{thm:blglob-cartier} gives the
following. Let $\cC_\cZ$ be the category whose objects are
$X$-schemes $T$ such that $f^*\cZ$ is Cartier, and whose morphisms are
morphisms over $X$. Then:

\begin{corollary}\label{cor:terminal}
$\Bl_\cZ X$ is a terminal object of $\cC_\cZ$.
\end{corollary}

No assumption is made on $X$ or $T$.

\begin{remark}
For $\iota$ a singleton, this is \cite[Tag 0806]{stacks-project}.
\end{remark}

\subsection{Formulation with mathlib ideal sheaves}\label{ssec:mathlib}

For an ideal sheaf $\cI$ on $X$, mathlib provides a closed subscheme
$V(\cI) \to X$ (\texttt{IdealSheafData.subscheme}); closed immersions into $X$
correspond to ideal sheaves on $X$ (\texttt{overEquivIdealSheafData}). For
$f : T \to X$, we write $f^{-1}\cI$ for the ideal sheaf of
$T \times_X V(\cI) \to T$ (\texttt{IdealSheafData.comap}).

\emph{From ideal sheaves to centres.} Let $(\cI_i)_{i \in \iota}$ be a family of
ideal sheaves on $X$. It defines a presented centre (\texttt{toPreClos}): the
schemes are the $V(\cI_i)$, the cover consists of the charts
$\Spec \sO_X(U) \to X$ for all affine opens $U$, the ideals are the $\cI_i(U)$,
and the isomorphisms of item (5) are \texttt{specQuotientIsoPullback}. We
write $[\cI] \in \Clos(X)$ for its class (\texttt{toClos}).

\emph{From centres to ideal sheaves.} Let $Z$ be a presented centre. Each
$Z_i \to X$ is a closed immersion: this is local on the target, and over the
chart $\Spec A_\gamma$ the base change $Z_i \times_X \Spec A_\gamma \to
\Spec A_\gamma$ is identified with $\Spec(A_\gamma/I_{i\gamma}) \to \Spec
A_\gamma$ (\texttt{PreClos.isClosedImmersion}). On $u_\gamma(U_\gamma)$, the
kernel of $Z_i \to X$ corresponds to $I_{i\gamma}$ under
$\sO_X(u_\gamma(U_\gamma)) \cong A_\gamma$ (\texttt{PreClos.ker\_ideal\_chart}).
Moreover, $[Z] = [\ker(Z_i \to X)]_i$ in $\Clos(X)$
(\texttt{PreClos.mk\_eq\_toClos\_ker}). A presented centre adds the choice of
an affine cover, on which the Cartier condition of Definition \ref{def:cartier}
is stated; $\Clos(X)$ forgets this choice.

\emph{Cartier conditions.} An ideal sheaf $\cI$ on $X$ is \emph{effective
Cartier} if every point of $X$ has an affine open neighbourhood $U$ such that
$\cI(U)$ is generated by a non-zero-divisor of $\sO_X(U)$.

\leansrc{Project/Blowups/IdealSheaf.lean}{202--203, 207--208}
\begin{lstlisting}
def IsPrincipalNZDOn (I : X.IdealSheafData) (U : X.affineOpens) : Prop :=
  ∃ f ∈ nonZeroDivisors Γ(X, U), I.ideal U = Ideal.span {f}

def IsEffectiveCartier (I : X.IdealSheafData) : Prop :=
  ∀ x : X, ∃ U : X.affineOpens, x ∈ U.1 ∧ I.IsPrincipalNZDOn U
\end{lstlisting}

\begin{proposition}\label{prop:cartier-mathlib}
Let $(\cI_i)_{i \in \iota}$ be a finite family of ideal sheaves on $X$ and
$f : T \to X$ a morphism. Then $f^*[\cI]$ is Cartier if and only if each
$f^{-1}\cI_i$ is effective Cartier (\texttt{isCars\_pullback\_toClos\_iff}).
\end{proposition}

\begin{proof}
Suppose each $f^{-1}\cI_i$ is effective Cartier. Generation by a
non-zero-divisor passes to basic opens, and two affine neighbourhoods of a
point contain a common basic open neighbourhood. Since $\iota$ is finite,
induction on $\iota$ shows that every point of $T$ has an affine open
neighbourhood on which all the $f^{-1}\cI_i$ are generated by non-zero-divisors
(\texttt{exists\_isPrincipalNZDOn}). Let $S$ be the set of such affine opens.
Presenting $(f^{-1}\cI_i)_i$ on the cover by $S$ gives a Cartier presented
centre (\texttt{toPreClosOn}, \texttt{isPreCars\_of\_generators}). It is
equivalent to $f^*[\cI]$, through the isomorphisms
$V(f^{-1}\cI_i) \cong T \times_X V(\cI_i)$ over $T$.

Conversely, let $W$ be a Cartier presented centre equivalent to $f^*[\cI]$, with
bijection $\sigma$ of index sets. The kernel of $W_k \to T$ is
$f^{-1}\cI_{\sigma(k)}$. For $x \in T$, take a chart of $W$ whose image
contains $x$. By \texttt{PreClos.ker\_ideal\_chart}, the sections of
$f^{-1}\cI_{\sigma(k)}$ on that image correspond to the ideal of $W$ on the
chart, which is generated by a non-zero-divisor.
\end{proof}

\begin{corollary}\label{cor:mathlib}
Let $(\cI_i)_{i \in \iota}$ be a finite family of ideal sheaves on $X$. There
is a scheme $B$ over $X$ such that each $\cI_i$ pulls back to an effective
Cartier ideal sheaf on $B$, and every $X$-scheme with this property has a
unique morphism to $B$ over $X$.
\end{corollary}

\leansrc{Project/Blowups/IdealSheaf.lean}{357--360}
\begin{lstlisting}
theorem exists_blowup_isEffectiveCartier [Finite ι] :
    ∃ (B : Scheme.{u+1}) (_ : B.Over X), (∀ i, ((I i).comap (B ↘ X)).IsEffectiveCartier) ∧
      ∀ (T : Scheme.{u+1}) [T.Over X], (∀ i, ((I i).comap (T ↘ X)).IsEffectiveCartier) →
        ∃! φ : T ⟶ B, φ.IsOver X := by
\end{lstlisting}

\begin{proof}
Apply Theorem \ref{thm:main} to $[\cI]$ and use Proposition
\ref{prop:cartier-mathlib}.
\end{proof}

For a finite family of closed immersions $Z_i \to X$, the same holds with
$\cI_i = \ker(Z_i \to X)$
(\texttt{exists\_blowup\_isEffectiveCartier\_of\_isClosedImmersion}).

\section{Dilatation}\label{sec:dilatation}

This section contains the construction that drives everything, following
\cite{mayeux-dilatations} and \cite{mayeux-zhang}.

\subsection{Multicenters}

\begin{definition}\label{def:multicenter}
Let $A$ be a commutative ring. A \emph{multicenter} on $A$ is a triple
$F = (\Lambda, (I_j)_{j \in \Lambda}, (a_j)_{j \in \Lambda})$ where $\Lambda$
is a type, each $I_j$ is an ideal of $A$ and each $a_j$ is an element of $A$.
The \emph{large ideal} at $j$ is
\[
  \cL_j \;:=\; I_j + (a_j).
\]
\end{definition}

\leansrc{Project/Dilatation/Multicenter.lean}{19--23, 39}
\begin{lstlisting}
@[ext]
structure Multicenter : Type (u+1) where
  (index : Type)
  (ideal : index → Ideal A)
  (elem : index → A)

def LargeIdeal (i : F.index) : Ideal A := F.ideal i + Ideal.span {F.elem i}
\end{lstlisting}

The constructions below depend on $F$ only through the pairs $(\cL_j, a_j)$;
$I_j$ enters only via $\cL_j$. Ideals $\cL_j$ with elements $a_j \in \cL_j$,
for $j \in \Lambda$, form a multicenter (take $I_j = \cL_j$), whose large
ideals are the $\cL_j$. The form $(I_j, a_j)$ is the one of
\cite{mayeux-dilatations} and of the formalisation \cite{mayeux-zhang}.

For $\nu \in \Lambda^{(\Nat)}$, a finitely supported function $\Lambda \to \Nat$,
write
\[
  \cL^\nu := \prod_{j} \cL_j^{\,\nu_j}, \qquad
  a^\nu := \prod_j a_j^{\,\nu_j},
\]
both finite products.

\subsection{The dilatation algebra}

\begin{definition}\label{def:dilatation}
The \emph{dilatation} of $A$ along $F$ is
\[
  A[F] \;:=\;
  \Bigl\{\, \tfrac{m}{a^\nu} \;:\; \nu \in \Lambda^{(\Nat)},\;
     m \in \cL^\nu \,\Bigr\} \Big/ \sim,
\]
where
\[
  \frac{m}{a^{\nu}} \sim \frac{n}{a^{\mu}}
  \iff
  \exists\, \beta \in \Lambda^{(\Nat)} : \quad
  m \, a^{\beta + \mu} \;=\; n \, a^{\beta + \nu} .
\]
\end{definition}

Addition and multiplication are the obvious ones,
\[
  \frac{m}{a^\nu} + \frac{n}{a^\mu} = \frac{m a^\mu + n a^\nu}{a^{\nu+\mu}},
  \qquad
  \frac{m}{a^\nu} \cdot \frac{n}{a^\mu} = \frac{mn}{a^{\nu+\mu}},
\]
and one checks these are well defined. Note that $m a^\mu \in
\cL^{\nu+\mu}$ because $a_j \in \cL_j$, which is where the definition of the
large ideal is used. The map $A \to A[F]$, $x \mapsto x/a^0$, makes $A[F]$ an
$A$-algebra.

\leansrc{Project/Dilatation/Multicenter.lean}{50--56, 89--90}
\begin{lstlisting}
structure PreDil where
  pow : F.index →₀ ℕ
  num : A
  num_mem : num ∈ F.LargeIdeal ^ pow

def r : F.PreDil → F.PreDil → Prop := fun x y =>
  ∃ β : F.index →₀ ℕ, x.num * F.elem^(β + y.pow) = y.num * F.elem^(β + x.pow)

variable (F) in
def Dilatation := Quotient F.setoid
\end{lstlisting}

Informally, $A[F] = A\bigl[\cL_j / a_j : j \in \Lambda\bigr]$, the subring of
the total localisation generated by the fractions $\ell/a_j$ with
$\ell \in \cL_j$. The quotient presentation above avoids having to embed
anything into a localisation.

\subsection{The two Cartier lemmas}

The following pair of statements is the reason dilatations are the right local
model. They say that in $A[F]$ the large ideals become principal, generated by
non-zero-divisors.

\begin{lemma}\label{lem:dil-nzd}
For every $\nu$, the image of $a^\nu$ in $A[F]$ is a non-zero-divisor.
\end{lemma}

\begin{lemma}\label{lem:dil-gen}
For every $\nu$, inside $A[F]$ one has
\[
  \bigl( a^{\nu} \bigr) \;=\; \cL^{\nu} \cdot A[F].
\]
\end{lemma}

\leansrc{Project/Dilatation/Multicenter.lean}{374, 386--388}
\begin{lstlisting}
lemma nonzerodiv_image (v : F.index →₀ ℕ) : algebraMap A A[F] (F.elem^v) ∈ nonZeroDivisors A[F] := by

lemma image_elem_LargeIdeal_equal (v : F.index →₀ ℕ) :
 Ideal.span ({algebraMap A A[F] (F.elem^v)}) =
    Ideal.map (algebraMap A A[F]) (F.LargeIdeal^v) := by
\end{lstlisting}

The proof of Lemma \ref{lem:dil-nzd} is short and shows what the relation is
doing. Suppose $\frac{m}{a^\mu} \cdot a^\nu = 0$ in $A[F]$. By definition of
$\sim$ there is $\beta$ with $m \, a^{\nu} a^{\beta} = 0$ in $A$. But then
taking $\beta' = \nu + \beta$ witnesses $\frac{m}{a^\mu} = 0$. So the element
we inverted is a non-zero-divisor not because it was one in $A$, which it need
not be, but because the equivalence relation was set up to kill exactly the
elements it annihilates.

Lemma \ref{lem:dil-gen} is the inclusion $\cL^\nu A[F] \subseteq (a^\nu)$,
which holds since $m \in \cL^\nu$ gives $m = a^\nu \cdot \frac{m}{a^\nu}$,
together with the reverse inclusion, which holds since $a^\nu \in \cL^\nu$.

\subsection{The universal property}\label{ssec:univ-dil}

\begin{theorem}\label{thm:dil-univ}
Let $B$ be an $A$-algebra such that
\begin{enumerate}
\item the image of $a_j$ in $B$ is a non-zero-divisor for every $j$, and
\item $(a_j) = \cL_j B$ inside $B$, for every $j$.
\end{enumerate}
Then there is exactly one morphism of $A$-algebras $A[F] \to B$.
\end{theorem}

\leansrc{Project/Dilatation/Multicenter.lean}{545--548, 622--625}
\begin{lstlisting}
def desc [Algebra A B]
    (non_zero_divisor : ∀ i : F.index, (algebraMap A B) (F.elem i) ∈ nonZeroDivisors B)
    (gen : ∀ i, Ideal.span {(algebraMap A B) (F.elem i)} = Ideal.map (algebraMap A B) (F.LargeIdeal i)) :
    A[F] →ₐ[A] B where

lemma lemma_exists_unique_morphism [Algebra A B]
    (non_zero_divisor : ∀ i : F.index, (algebraMap A B) (F.elem i) ∈ nonZeroDivisors B)
    (gen : ∀ i, Ideal.span {(algebraMap A B) (F.elem i)} = Ideal.map (algebraMap A B) (F.LargeIdeal i))
    (χ' : A[F] →ₐ[A] B) : χ' = desc F non_zero_divisor gen := by
\end{lstlisting}

Existence and uniqueness both come from the same observation. Given $\nu$ and
$m \in \cL^\nu$, hypothesis (2) gives $m \in \cL^\nu B = (a^\nu)$, so there is
$b \in B$ with $m = a^\nu b$, and hypothesis (1) says this $b$ is unique. Send
$\frac{m}{a^\nu}$ to that $b$. Any $A$-algebra map $\chi$ must satisfy
$a^\nu \chi\bigl(\frac{m}{a^\nu}\bigr) = \chi(m) = m$, so it must send
$\frac{m}{a^\nu}$ to the same $b$. Existence and uniqueness are the two halves
of ``divide by $a^\nu$, and the quotient is unique''.

The hypotheses on $B$ are the Cartier condition with prescribed generators.
So Theorem \ref{thm:dil-univ} reads:

\begin{quote}
$A[F]$ is the initial $A$-algebra $B$ such that, for every $j$,
$\cL_j B = (a_j)$ with $a_j$ a non-zero-divisor in $B$.
\end{quote}

That is the affine case of Corollary \ref{cor:terminal}, in ring-theoretic
form. Everything in Section \ref{sec:proof} is an elaboration of this one
statement: first upgrading it from a single multicenter to a family of ideals
with all of its presentations, then from affine to global.

The following special case follows from Theorem \ref{thm:dil-univ} by
uniqueness.

\begin{corollary}\label{cor:dil-trivial}
If every $a_j$ is already a non-zero-divisor in $A$ and $\cL_j = (a_j)$, then
$A \to A[F]$ is an isomorphism.
\end{corollary}

Indeed $B = A$ then satisfies the hypotheses of Theorem \ref{thm:dil-univ}, so
there is a unique $A[F] \to A$, and both composites with $A \to A[F]$ are
identities by the uniqueness clause applied to $A$ and to $A[F]$ in turn.

\leansrc{Project/Dilatation/Multicenter.lean}{657--658}
\begin{lstlisting}
def ofFamilyIso {index : Type} {c : index → A} (c0 : ∀ i, c i ∈ nonZeroDivisors A) :
    A[ofFamily index c] ≃ₐ[A] A :=
\end{lstlisting}

\subsection{Functoriality}\label{ssec:dil-funct}

Let $A \to B$ be a ring map and $F$ a multicenter on $A$. Pushing forward the
ideals and the elements gives a multicenter $F_B$ on $B$ with the same index
type, ideals $I_j B$ and elements the images of $a_j$. One checks
$\cL_j(F_B) = \cL_j(F) B$, so the large ideals also push forward.

\begin{proposition}\label{prop:dil-funct}
There is a unique morphism of $A$-algebras $A[F] \to B[F_B]$ compatible with
$A \to B$.
\end{proposition}

Again this is Theorem \ref{thm:dil-univ}: in $B[F_B]$ the element $a_j$ is a
non-zero-divisor generating $\cL_j(F) B[F_B]$, by Lemmas \ref{lem:dil-nzd} and
\ref{lem:dil-gen} applied over $B$. So $A[F] \to B[F_B]$ exists and is unique.
The pushed forward multicenter $F_B$ is \texttt{image\_mult}, and the map is
built by \texttt{desc}, that is by the universal property.

\leansrc{Project/Dilatation/Multicenter.lean}{772--775, 784}
\begin{lstlisting}
def image_mult [Algebra A B] : Multicenter B :=
  { index := _
    ideal i := Ideal.map (algebraMap A B) (F.ideal i)
    elem i := algebraMap A B (F.elem i)}

def functo_dila_alg [Algebra A B] : A[F] →ₐ[A] B[image_mult (B := B) F] :=
  desc F ...
\end{lstlisting}

When $A \to B$ is flat one can say more: the natural map
\[
  B \otimes_A A[F] \longrightarrow B[F_B]
\]
is an isomorphism. Flatness is genuinely needed here, since the construction of
the inverse uses that a non-zero-divisor stays a non-zero-divisor after a flat
base change.

\subsection{Presentations}

We now fix a finite index type $\iota$ and a family $L = (L_i)_{i \in \iota}$
of ideals of $A$. This is the affine centre we want to blow up. A single
multicenter is not enough, because a multicenter comes with a choice of the
elements $a_j$, and there is no canonical choice. The construction takes all
choices at once.

\begin{definition}\label{def:pres}
A \emph{presentation} of $L$ is a multicenter $F$ on $A$ with finite index type
$\Lambda$, together with maps $\Psi : \Lambda \to \iota$ and
$s : \iota \to \Lambda$ such that
\[
  \Psi \circ s = \mathrm{id}_\iota,
  \qquad
  \cL_j = L_{\Psi(j)} \quad \text{for all } j \in \Lambda .
\]
We write $\Pres(L)$ for the collection of presentations.
\end{definition}

\leansrc{Project/Blowups/Bl.lean}{18--24}
\begin{lstlisting}
structure Mu : Type (u + 1) where
multicenter : Multicenter A
[fin : Fintype multicenter.index]
Ψ : multicenter.index → ι
sec : ι → multicenter.index
surj : ∀ i, Ψ (sec i) = i
cond : ∀ i, multicenter.LargeIdeal i = L (Ψ i)
\end{lstlisting}

Unwinding, a presentation is a finite set $\Lambda$ of ``chart directions'',
each labelled by the ideal $L_{\Psi(j)}$ it refers to, and for each $j$ an
element $a_j \in L_{\Psi(j)}$, presented as
\[
  L_{\Psi(j)} \;=\; I_j + (a_j).
\]
The section $s$ exists if and only if $\Psi$ is surjective. Nothing
requires $\Psi$ to be injective, so an ideal may appear with several different
chosen elements.

\subsection{Dilatations as charts of a multigraded Proj}\label{ssec:proj}

The dilatations $A[F]$, $F \in \Pres(L)$, are the charts of a single scheme.
To glue them we realise each as an affine piece of a multigraded $\Proj$.

Let
\[
  \Rees(L) \;:=\; \bigoplus_{v \in \iota^{(\Nat)}} L^v ,
  \qquad
  L^v := \prod_i L_i^{\,v_i},
\]
the multi-Rees algebra of the family, graded by $\Nat^{\iota}$ and then by
$\Int^\iota$ through the inclusion $\Nat^\iota \hookrightarrow \Int^\iota$. For
a homogeneous submonoid $S$ of $\Rees(L)$, let
\[
  \Pot(S) \;:=\; \Bigl( S^{-1} \Rees(L) \Bigr)_0
\]
be the degree zero part of the homogeneous localisation, called the
\emph{potion} of $S$ in the formalisation. The multigraded $\Proj$ is glued
from the affine schemes $\Spec \Pot(S)$ as $S$ runs over the relevant finitely
generated homogeneous submonoids.
Here the \emph{saturation} $\bar S$ of $S$ is the submonoid of homogeneous
elements of $\Rees(L)$ dividing some element of $S$, and $S$ is
\emph{relevant} if for every $v \in \Int^\iota$ some positive multiple $nv$
lies in the subgroup generated by the degrees of $\bar S$.

Given a presentation $F \in \Pres(L)$, let $S_F \subseteq \Rees(L)$ be the
submonoid generated by the elements
\[
  a_j \in L_{\Psi(j)} = \bigl(\Rees L\bigr)_{e_{\Psi(j)}},
  \qquad j \in \Lambda,
\]
that is, $a_j$ placed in degree $e_{\Psi(j)}$, the $\Psi(j)$-th basis vector.

\begin{proposition}\label{prop:dil-is-potion}
There is a canonical isomorphism of $A$-algebras
\[
  A[F] \;\xrightarrow{\ \sim\ }\; \Pot(S_F).
\]
\end{proposition}

A homogeneous element of $S_F$ of degree $\nu \in \Nat^\iota$ is a
product $a^\lambda$ for $\lambda \in \Lambda^{(\Nat)}$ with
$\sum_j \lambda_j e_{\Psi(j)} = \nu$. A degree zero element of the localisation
$S_F^{-1}\Rees(L)$ is therefore a fraction
\[
  \frac{m}{a^{\lambda}}, \qquad m \in \Rees(L)_\nu = L^{\nu},
\]
and the defining condition $\cL_j = L_{\Psi(j)}$ of a presentation says exactly
that
\[
  L^{\nu} \;=\; \prod_i L_i^{\nu_i} \;=\; \prod_j \cL_j^{\lambda_j}
  \;=\; \cL^{\lambda}.
\]
So the numerators available on the two sides coincide. The equivalence relation
also matches: two fractions in a homogeneous localisation agree when they agree
after multiplying by some element of $S_F$, and elements of $S_F$ are exactly
the $a^\beta$, which is the relation of Definition \ref{def:dilatation}. The
map itself is again given by the universal property of the dilatation, and the
isomorphism is obtained from injectivity and surjectivity, proved separately.

\leansrc{Project/Blowups/Bl.lean}{114--116, 271--272, 605--607}
\begin{lstlisting}
def clo_mu (P : Mu L) :
    HomogeneousSubmonoid (ReesAlgebra.intGrading L) :=
  HomogeneousSubmonoid.closure

def clo_mu_mor (P : Mu L) : A[P.multicenter] →ₐ[A] (clo_mu L P).Potion :=
  Multicenter.desc P.multicenter

def Mu_mor_iso (P : Mu L) :
    A[P.multicenter] ≃ₐ[A] (clo_mu L P).Potion :=
  AlgEquiv.ofBijective (clo_mu_mor ..) <| ⟨clo_mu_mor_inj L P, clo_mu_mor_surj L P⟩
\end{lstlisting}

The inclusion $\cL_j \subseteq L_{\Psi(j)}$ is used to define the map, and the
reverse inclusion for surjectivity.

For the gluing we also need to know that $S_F$ is an admissible index for the
$\Proj$ construction, meaning finitely generated and relevant. Finite
generation is clear, since $\Lambda$ is finite and $S_F$ is generated by the
$a_j$. Relevance amounts to the degrees of $\bar S_F$ generating a finite-index
subgroup of $\Int^\iota$; they generate $\Int^\iota$, since $S_F \subseteq \bar S_F$,
$a_j$ has degree $e_{\Psi(j)}$, and $\Psi$ is surjective.

\begin{definition}\label{def:blm}
$\BlM(L) := \Proj$ of $\Rees(L)$ taken over the family
$\{ S_F : F \in \Pres(L)\}$, that is, the scheme obtained by gluing the affine
charts $\Spec \Pot(S_F) \cong \Spec A[F]$, $F \in \Pres(L)$, along the standard
$\Proj$ gluing data. It carries a structure morphism $\BlM(L) \to \Spec A$
induced by $A \cong \Rees(L)_0$.
\end{definition}

\leansrc{Project/Blowups/Bl.lean}{703--706, 721--722, 726--728}
\begin{lstlisting}
def map_index (P : Mu L) :
    GoodPotionIngredient (ReesAlgebra.intGrading L) where
  toHomogeneousSubmonoid := clo_mu L P
  relevant := clo_mu_rel L P

def BlMu : Scheme :=
  Proj (τ := Mu L) (map_index L)

instance BlMuOverSpec : Scheme.Over (BlMu L) (Spec <| CommRingCat.of A) where
  hom := (GoodPotionIngredient.over (ℱ := map_index L)) ≫
    Spec.map (CommRingCat.ofHom <| ReesAlgebra.degreeZeroIso' L)
\end{lstlisting}

The field \texttt{relevant} is the relevance condition discussed above, and the
structure morphism to $\Spec A$ is induced by the isomorphism
$A \cong \Rees(L)_0$ recorded as \texttt{degreeZeroIso'}.

Two remarks. First, the gluing data is inherited from the ambient multigraded
$\Proj$, so the overlaps and the cocycle condition are inherited from the
$\Proj$ construction, requiring no additional verification. Second, we do not need to know
that $\BlM(L)$ is all of $\Proj \Rees(L)$. It is the union of the dilatation
charts, and every statement we prove is scoped to those charts, reflecting a
deliberate architectural choice that avoids the need for a comparison that would
otherwise complicate the proof.

To summarise the picture:
\[
\begin{tikzcd}[column sep=large]
  \Spec A[F] \arrow[r,"\sim"] \arrow[drr, bend right=15] & \Spec \Pot(S_F) \arrow[r,hook] & \BlM(L) \arrow[d] \\
   & & \Spec A
\end{tikzcd}
\]
Each presentation gives a dilatation, each dilatation is an affine chart, and
the charts glue to the affine blowup.

\section{Proof of the theorem}\label{sec:proof}

The proof has two halves. First the affine case: $\BlM(L)$ is terminal among
affine-based $\Spec A$-schemes on which $L$ becomes Cartier. Then the global
case: the local blowups over the charts of $X$ are glued, and the gluing is
governed at every stage by the affine uniqueness statement.

\subsection{The affine centre is Cartier on \texorpdfstring{$\BlM(L)$}{Bl(L)}}
\label{ssec:cars}

The first thing to prove is that $\BlM(L)$ does what it is supposed to do.

\begin{proposition}\label{prop:blm-cartier}
The centre $L$ pulled back to $\BlM(L)$ is Cartier.
\end{proposition}

The presentation exhibiting this is: take as cover the family
of dilatation charts $\Spec A[F]$, $F \in \Pres(L)$, indexed by $\Pres(L)$
itself, and as ideal on the chart $F$ the image $L_i \cdot A[F]$.

For each $i$ and each $F$, we must produce a generator of $L_i \cdot A[F]$ that
is a non-zero-divisor. Choose $j$ with $\Psi(j) = i$ (in the formalisation,
$j = s(i)$), so $\cL_j = L_i$. Then Lemmas
\ref{lem:dil-nzd} and \ref{lem:dil-gen}, applied with $\nu = e_j$, give
\[
  L_i \cdot A[F] \;=\; \cL_j \cdot A[F] \;=\; \bigl( a_j \bigr),
  \qquad a_j \text{ a non-zero-divisor in } A[F].
\]
So the generator is $a_j$. The formal proof passes these two lemmas to a
criterion that checks the Cartier property chart by chart.

\leansrc{Project/Blowups/Cars.lean}{64--68}
\begin{lstlisting}
lemma BlMuPreClos_IsPreCars (A : CommRingCat) (L : ι → Ideal A) :
    IsPreCars _ (BlMuPreClos A L) := by
  refine isPreCars_of_generators (BlMuPreClos A L)
      (fun i (P : Mu L) => (Mu_mor_iso L P).toRingEquiv
        (algebraMap A A[P.multicenter] (P.multicenter.elem (P.sec i)))) ?_ ?_
\end{lstlisting}

The presented centre and the Cartier statement are as follows. Its cover is \texttt{BlMuAffineCover}, the family of dilatation
charts, and its ideal on the chart $P$ is the image of $L_i$.

\leansrc{Project/Blowups/Cars.lean}{43--48, 89--91}
\begin{lstlisting}
def BlMuPreClos (A : CommRingCat) (L : ι → Ideal A) : PreClos (BlMu L) where
  indnumb := ι
  subscheme i := pullback ((loc_to_PreClos A L).subscheme i ↘ Spec A) (BlMu L ↘ Spec A)
  over i := ⟨pullback.snd _ _⟩
  cov := BlMuAffineCover A L
  ideal i P := Ideal.map (algebraMap A ((map_index L P).Potion)) (L i)

lemma blowups_Cars (A : CommRingCat) (L : ι → Ideal A) :
    IsCars _ <|
      pullback_Clos (BlMu L ↘ (Spec (CommRingCat.of A))) (loc_to_Clos A L) :=
\end{lstlisting}

The term \texttt{P.multicenter.elem (P.sec i)} is $a_j$ for $j = s(i)$, and
\texttt{Mu\_mor\_iso} is the isomorphism of Proposition \ref{prop:dil-is-potion}
transporting it into the potion.

A consequence we use repeatedly: an open subscheme of $\BlM(L)$ still has the
property, since open immersions are flat and the Cartier property is stable under flat
base change. In the formalisation this is
\texttt{open\_of\_blowup\_IsCars} (\texttt{Project/Blowups/Cars.lean}, lines
100--107).

\subsection{Affine uniqueness}\label{ssec:aff-uniq}

\begin{proposition}\label{prop:aff-uniq}
Let $T$ be a scheme over $\Spec A$ on which $L$ becomes Cartier. Then any two
morphisms $T \to \BlM(L)$ over $\Spec A$ are equal.
\end{proposition}

\begin{proof}
Let $p : T \to \Spec A$ be the structure morphism, let
$j_F : \Spec A[F] \cong \Spec \Pot(S_F) \to \BlM(L)$, $F \in \Pres(L)$, be the
chart inclusions of Definition \ref{def:blm}, and let
$\varphi, \varphi' : T \to \BlM(L)$ be morphisms over $\Spec A$. Two morphisms
of schemes that agree on an open neighbourhood of every point are equal, so it
suffices to fix $x \in T$ and find such a neighbourhood.

By hypothesis and Definition \ref{def:cartier}, the pulled back centre has a
Cartier presentation $Z$, with an affine open cover
$(u_\gamma : \Spec B_\gamma \to T)_\gamma$, and by Definition \ref{def:rel}
there is a bijection $\sigma$ from the index type of $Z$ to $\iota$. Since
$p \circ u_\gamma = \Spec(A \to B_\gamma)$, the ideal of $Z$ of index $k$ on
$B_\gamma$ equals $L_{\sigma(k)} B_\gamma$
(\texttt{PreClos\_ideal\_eq\_of\_rel}). Reindexing by $\sigma$, for each
$i \in \iota$ and each $\gamma$ there is a non-zero-divisor
$z_{i\gamma} \in B_\gamma$ with $L_i B_\gamma = (z_{i\gamma})$.

Choose $\gamma$ with $x \in u_\gamma(\Spec B_\gamma)$, and
$F, G \in \Pres(L)$ with $\varphi(x) \in j_F(\Spec A[F])$ and
$\varphi'(x) \in j_G(\Spec A[G])$. Choose an affine open $\Spec B \ni x$
contained in
$\varphi^{-1}(j_F(\Spec A[F])) \cap \varphi'^{-1}(j_G(\Spec A[G])) \cap u_\gamma(\Spec B_\gamma)$.
Since $j_F$, $j_G$, $u_\gamma$ are open immersions, we obtain ring maps
$g : A[F] \to B$, $g' : A[G] \to B$, $h : B_\gamma \to B$ with
$\varphi|_{\Spec B} = j_F \circ \Spec(g)$ and
$\varphi'|_{\Spec B} = j_G \circ \Spec(g')$. As all morphisms are over
$\Spec A$, $g$ and $g'$ are $A$-algebra maps and $A \to B$ factors through $h$.

The map $h$ is flat, because $\Spec B \to \Spec B_\gamma$ is an open immersion
of affine schemes (Section \ref{ssec:cars}). Hence $c_i := h(z_{i\gamma})$ is
a non-zero-divisor of $B$ and $L_i B = (c_i)$ for every $i \in \iota$.

Let $F \sqcup G \in \Pres(L)$ be the presentation indexed by
$\Lambda_F \sqcup \Lambda_G$, carrying the pairs of $F$ and of $G$, with
section that of $F$ (\texttt{union\_Mu}), and let
$\lambda_F : A[F] \to A[F \sqcup G]$ and $\lambda_G : A[G] \to A[F \sqcup G]$
be the structure maps (\texttt{dilationToUnion\_left},
\texttt{dilationToUnion\_right}). For $j \in \Lambda_F$, applying $g$ to Lemma \ref{lem:dil-gen} gives
$(a_j) = \cL_j B = L_{\Psi(j)} B = (c_{\Psi(j)})$ in $B$. Hence $a_j$ and the
non-zero-divisor $c_{\Psi(j)}$ generate the same principal ideal, so $a_j$ is a
non-zero-divisor of $B$. The same holds for $j \in \Lambda_G$, using $g'$. By
Theorem \ref{thm:dil-univ} there is an $A$-algebra map
$g'' : A[F \sqcup G] \to B$. The maps $g$ and $g'' \circ \lambda_F$ are
$A$-algebra maps $A[F] \to B$, so they are equal by the uniqueness in Theorem
\ref{thm:dil-univ}; likewise $g' = g'' \circ \lambda_G$
(\texttt{lemm\_dila\_double\_union}).

One has $S_F, S_G \le S_{F \sqcup G}$, and under the isomorphisms of
Proposition \ref{prop:dil-is-potion} the localisation maps
$\Pot(S_F) \to \Pot(S_{F \sqcup G})$ and $\Pot(S_G) \to \Pot(S_{F \sqcup G})$
are $\lambda_F$ and $\lambda_G$ (\texttt{Mu\_mor\_iso\_commutes\_ringHom'},
\texttt{Mu\_mor\_iso\_commutes\_ringHom\_right'}). The gluing data of
$\BlM(L)$ then give
$j_F \circ \Spec(\lambda_F) = j_{F \sqcup G} = j_G \circ \Spec(\lambda_G)$
(\texttt{proj\_glue\_condition}). Therefore
\[
  \varphi|_{\Spec B}
  = j_F \circ \Spec(\lambda_F) \circ \Spec(g'')
  = j_{F \sqcup G} \circ \Spec(g'')
  = j_G \circ \Spec(\lambda_G) \circ \Spec(g'')
  = \varphi'|_{\Spec B}.
\]
\end{proof}

\leansrc{Project/Blowups/UniqueBlowup.lean}{254--261}
\begin{lstlisting}
theorem ProjBlowup_UnivProp_unicity_affine
  (A : CommRingCat.{u+1}) (L : ι → Ideal A) [fin : Fintype ι]
  {T : Scheme} [T.Over (Spec A)]
  (cond : IsCars _ (pullback_Clos (T ↘ Spec A) (loc_to_Clos A L)))

  (φ φ' : T ⟶ BlMu L)
  (φ_over : Scheme.Hom.IsOver φ (Spec A))
  (φ'_over : Scheme.Hom.IsOver φ' (Spec A)) : φ = φ' := by
\end{lstlisting}

\subsection{Affine existence}\label{ssec:aff-exist}

\begin{proposition}\label{prop:aff-exist}
Let $T$ be a scheme over $\Spec A$ on which $L$ becomes Cartier. Then there is
a morphism $T \to \BlM(L)$ over $\Spec A$.
\end{proposition}

\begin{proof}
We follow \texttt{ProjBlowup\_UnivProp\_existence\_affine\_preclo}
(\texttt{Project/Blowups/BlowupExists.lean}). Write $e_i\in\Int^\iota$ for the
basis vectors and $d=\sum_i e_i$. For a homogeneous element $f$ of a Rees
algebra, $\langle f\rangle$ denotes the submonoid it generates. Saturation and relevance are as in Section \ref{ssec:proj}; a homogeneous submonoid and its saturation have the same potion (\texttt{equivBarPotion}). If $S$ is generated by finitely many
homogeneous elements with product $f$, then $\bar S=\overline{\langle f\rangle}$,
since each generator divides $f$ and $f\in S$; hence $\Pot(S)=\Pot\langle f\rangle$.

\emph{Step 1 (local data).} Choose a Cartier presentation of the pullback of $L$
to $T$ (Definition \ref{def:cartier}). It gives an affine open cover
$(u_\gamma:\Spec B_\gamma\to T)_\gamma$ and non-zero-divisors $c_i\in B_\gamma$
with $L_iB_\gamma=(c_i)$, after reindexing by the bijection of Definition \ref{def:rel} (\texttt{PreClos\_ideal\_eq\_of\_rel}). Fix $\gamma$ and
write $B=B_\gamma$; the map $\Spec B\to\Spec A$ is $\Spec$ of a ring map $A\to B$.

\emph{Step 2 ($B$ as a potion).} The pairs $(L_iB,c_i)$, with
$\Psi=s=\mathrm{id}$, form a presentation $F_B$ of $LB$ over $B$, since
$\cL_i=L_iB+(c_i)=(c_i)=L_iB$. Put $c=\prod_i c_it^{e_i}\in\Rees(LB)$, of degree
$d$. By Corollary \ref{cor:dil-trivial}, Proposition \ref{prop:dil-is-potion} and
the preliminary remark, $B\cong B[F_B]\cong\Pot(S_{F_B})=\Pot\langle c\rangle$ as
$B$-algebras.

\emph{Step 3 (expansion).} Since $c_i\in L_iB$, there are a finite set
$t_i\subseteq L_i$ and $\lambda_i(a)\in B$ with $c_i=\sum_{a\in t_i}\lambda_i(a)a$
in $B$. For $p\in P=\prod_i t_i$ put $g_p=\prod_i p_it^{e_i}\in\Rees(LB)$ and
$f_p=\prod_i\lambda_i(p_i)p_it^{e_i}$. Then $c=\sum_{p\in P}f_p$, and each $f_p$
is homogeneous of degree $d$. It is relevant, since each
$\lambda_i(p_i)p_it^{e_i}$ divides $f_p$ and the $e_i$ generate $\Int^\iota$. The
elements $f_p/c$ of $\Pot\langle c\rangle$ sum to $1$, so the opens
$D(f_p/c)\cong\Spec\Pot\langle cf_p\rangle$ cover $\Spec\Pot\langle c\rangle$.
Each maps to $\Spec\Pot\langle f_p\rangle$, and these maps glue to a morphism over
$\Spec B$ from $\Spec\Pot\langle c\rangle$ to the scheme $U$ glued from the
$\Spec\Pot\langle f_p\rangle$, $p\in P$, along the $\Proj$ gluing data of
$\Rees(LB)$ (\texttt{sum\_open\_finset}).

\emph{Step 4 (charts of $\BlM(L)$).} Since $g_p$ divides $f_p$,
$\overline{\langle g_p\rangle}\subseteq\overline{\langle f_p\rangle}$, which gives
a map $\Pot\langle g_p\rangle\to\Pot\langle f_p\rangle$. Let
$\tilde g_p=\prod_ip_it^{e_i}\in\Rees(L)$; the graded map $\Rees(L)\to\Rees(LB)$
sends $\tilde g_p$ to $g_p$ and gives
$\Pot\langle\tilde g_p\rangle\to\Pot\langle g_p\rangle$. The pairs $(L_i,p_i)$,
with $\Psi=s=\mathrm{id}$, form a presentation $F_p\in\Pres(L)$ (Definition
\ref{def:pres}), since $\cL_i=L_i+(p_i)=L_i=L_{\Psi(i)}$ as $p_i\in L_i$. By the
preliminary remark, $\Pot(S_{F_p})=\Pot\langle\tilde g_p\rangle$. With the chart
inclusion $j_{F_p}$ of Definition \ref{def:blm} we obtain
\[
h_p:\Spec\Pot\langle f_p\rangle\to\Spec\Pot\langle g_p\rangle\to
\Spec\Pot(S_{F_p})\xrightarrow{j_{F_p}}\BlM(L),
\]
a morphism over $\Spec A$.

\emph{Step 5 (gluing on $U$).} Let $\Spec\Pot(S)$ be the overlap of two charts
of $U$. The image of $c_i$ in $\Pot(S)$ is a non-zero-divisor
(\texttt{ReesAlgebra.potion\_algebraMap\_mem\_nonZeroDivisors}) and generates
$L_i\Pot(S)$, so $L$ becomes Cartier on $\Spec\Pot(S)$
(\texttt{loc\_isCars\_of\_principal}). The restrictions of $h_p$ and $h_q$ to the
overlap are both over $\Spec A$, so they agree by Proposition
\ref{prop:aff-uniq}. Hence the $h_p$ glue to $h:U\to\BlM(L)$ over $\Spec A$. Put
$\varphi_\gamma:\Spec B\cong\Spec\Pot\langle c\rangle\to U\xrightarrow{h}\BlM(L)$;
it is over $\Spec A$.

\emph{Step 6 (gluing on $T$).} For $\gamma,\delta$, the map
$W=\Spec B_\gamma\times_T\Spec B_\delta\to T$ is an open immersion, hence flat
(Section \ref{ssec:cars}), so $L$ becomes Cartier on $W$
(\texttt{pullback\_IsCars}). By Proposition \ref{prop:aff-uniq},
$\varphi_\gamma$ and $\varphi_\delta$ agree on $W$. Hence the $\varphi_\gamma$
glue to $\varphi:T\to\BlM(L)$ with $\varphi\circ u_\gamma=\varphi_\gamma$; it is
over $\Spec A$.
\end{proof}

Putting Propositions \ref{prop:blm-cartier}, \ref{prop:aff-uniq} and
\ref{prop:aff-exist} together:

\leansrc{Project/Blowups/BlowupExists.lean}{989--993, 1014--1018}
\begin{lstlisting}
lemma ProjBlowup_UnivProp_existence_affine
    (A : CommRingCat) (L : ι → Ideal A) [fin : Fintype ι]
    {T : Scheme} [T.Over (Spec A)]
    (cond : pullback_Clos (T ↘ Spec A) (loc_to_Clos A L) ∈ CarsAsSubsetOfClos T) :
    ∃ φ : T ⟶ BlMu L, Scheme.Hom.IsOver φ (Spec A) := by

lemma ProjBlowup_UnivProp_affine
  (A : CommRingCat) (L : ι → Ideal A) [fin : Fintype ι]
  {T : Scheme} [T.Over (Spec A)]
  (cond : pullback_Clos (T ↘ Spec A) (loc_to_Clos A L) ∈ CarsAsSubsetOfClos T) :
    ∃! φ : T ⟶ BlMu L, Scheme.Hom.IsOver φ (Spec A) := by
\end{lstlisting}

\begin{theorem}\label{thm:aff-main}
$\BlM(L)$ is terminal among $\Spec A$-schemes on which $L$ becomes Cartier.
\end{theorem}

The formalisation packages this as a structure, so that the affine blowup is an
object carrying its own universal property rather than a scheme with three
separate theorems attached.

\leansrc{Project/Blowups/UniqueBlowup.lean}{38--47}
\begin{lstlisting}
structure conceptual_blowup (Z : Clos X) where
  scheme : Scheme
  over : Scheme.Over scheme X
  in_cars : pullback_Clos (scheme ↘ X) Z ∈ CarsAsSubsetOfClos scheme
  φ (T : Scheme) [T.Over X] (in_preCars : pullback_Clos (T ↘ X) Z ∈ (CarsAsSubsetOfClos T)) :
    T ⟶ scheme
  φ_over (T : Scheme) [T.Over X] (in_preCars : pullback_Clos (T ↘ X) Z ∈ (CarsAsSubsetOfClos T)) :
    Scheme.Hom.IsOver (φ T in_preCars) X
  φ_uniq (T : Scheme) [T.Over X] (in_preCars : pullback_Clos (T ↘ X) Z ∈ (CarsAsSubsetOfClos T)) :
    ∀ φ' : T ⟶ scheme, Scheme.Hom.IsOver φ' X → φ' = φ T in_preCars
\end{lstlisting}

Theorem \ref{thm:aff-main} is then the statement that $\BlM(L)$ carries this
structure, for the centre $L$ viewed as an element of $\Clos(\Spec A)$.

\leansrc{Project/Blowups/BlowupExists.lean}{1042--1049}
\begin{lstlisting}
def ProjBlowup_is_conceptual_blowups_affine (A : CommRingCat) (L : ι → Ideal A) [fin : Fintype ι] :
    conceptual_blowup (loc_to_Clos A L) where
  scheme := BlMu L
  over := inferInstance
  in_cars := blowups_Cars A L
  φ T _ cond := ProjBlowup_φ A L cond
  φ_over T _ cond := ProjBlowup_φ_over A L cond
  φ_uniq T _ cond φ' hφ' := ProjBlowup_φ_uniq A L cond φ' hφ'
\end{lstlisting}

\subsection{Base change along open immersions}\label{ssec:bc}

To globalise we need to know how $\BlM$ behaves under restriction to an open.

\begin{proposition}\label{prop:bc}
Let $A \to B$ be a ring map with $\Spec B \to \Spec A$ an open immersion, and
let $L' := (L_i B)_i$. Then there is exactly one isomorphism
\[
  \BlM(L) \times_{\Spec A} \Spec B \;\cong\; \BlM(L')
\]
over $\Spec B$.
\end{proposition}

\leansrc{Project/Blowups/BlowupExists.lean}{1337--1341}
\begin{lstlisting}
lemma base_change_Bl_open [Fintype ι] (A B : CommRingCat) [Algebra A B]
    [IsOpenImmersion (Spec B ↘ Spec A)] (L : ι → Ideal A) :
    ∃! (e : pullback (BlMu L ↘ Spec A) (Spec B ↘ Spec A) ≅

    BlMu (L := fun i : ι => Ideal.map (algebraMap A B) (L i))), Scheme.Hom.IsOver e.hom (Spec B) := by
\end{lstlisting}

\begin{proof}
Put $P := \BlM(L) \times_{\Spec A} \Spec B$, with projections $\mathrm{pr}_1$,
$\mathrm{pr}_2$. To show that $L$ is Cartier on $P$, we apply a more general
result, namely that the Cartier property is stable under flat base change
(Section \ref{ssec:cars}; \texttt{pullback\_IsCars}): $\mathrm{pr}_1$ is a
base change of the open immersion $\Spec B \to \Spec A$, hence flat, and $L$
is Cartier on $\BlM(L)$ (Proposition \ref{prop:blm-cartier}). The pullback of
$L$ to $\Spec B$ is $L'$ (\texttt{loc\_to\_Clos\_baseChange}). Hence $L'$ is
Cartier on $P$, and, since $L'$ is Cartier on $\BlM(L')$ (Proposition
\ref{prop:blm-cartier}), $L$ is Cartier on $\BlM(L')$. Theorem
\ref{thm:aff-main} gives $\varphi : P \to \BlM(L')$ over $\Spec B$ and
$\psi : \BlM(L') \to \BlM(L)$ over $\Spec A$. The morphisms $\psi$ and
$\BlM(L') \to \Spec B$ have the same composite to $\Spec A$, hence define
$\psi' : \BlM(L') \to P$ over $\Spec B$. Then
$\mathrm{pr}_1 \circ \psi' \circ \varphi = \mathrm{pr}_1$ by the uniqueness in
Theorem \ref{thm:aff-main} for $L$ on $P$ over $\Spec A$, and
$\mathrm{pr}_2 \circ \psi' \circ \varphi = \mathrm{pr}_2$ since $\varphi$ is
over $\Spec B$; hence $\psi' \circ \varphi = \mathrm{id}_P$. By the uniqueness
for $L'$ on $\BlM(L')$ over $\Spec B$, $\varphi \circ \psi' = \mathrm{id}$. By
the uniqueness for $L'$ on $P$ over $\Spec B$, any isomorphism
$P \to \BlM(L')$ over $\Spec B$ equals $\varphi$.
\end{proof}

\subsection{The local blowups and their overlaps}\label{ssec:overlaps}

Now let $Z$ be a presented centre on $X$, with cover
$\cU = (u_\gamma : U_\gamma = \Spec A_\gamma \to X)_{\gamma \in J}$ and ideals
$I_{i\gamma}$. For each $\gamma$ set
\[
  B_\gamma \;:=\; \BlM\bigl( (I_{i\gamma})_{i} \bigr),
\]
the affine blowup of the chart, a scheme over $\Spec A_\gamma$ and hence over
$X$.

For a pair $\gamma, \delta$ let $U_{\gamma\delta} := U_\gamma \times_X U_\delta$
be the overlap, an open of both charts, and let
\[
  V_{\gamma\delta} \;:=\; B_\gamma \times_{\Spec A_\gamma} U_{\gamma\delta}
\]
be its preimage in $B_\gamma$. This is an open subscheme of $B_\gamma$, and we
write $f_{\gamma\delta} : V_{\gamma\delta} \hookrightarrow B_\gamma$ for the
inclusion and $\pi_{\gamma\delta} : V_{\gamma\delta} \to U_{\gamma\delta}$ for
the projection.

\leansrc{Project/Blowups/BlowupExists.lean}{1393--1400, 1417--1418, 1437--1441, 1614--1620}
\begin{lstlisting}
def ideal_loc (X : Scheme) (Z : PreClos X) (γ : Z.cov.J) : Z.indnumb → Ideal (Z.cov.obj γ) :=
  fun (i : Z.indnumb) => Z.ideal i γ

def Proj_loc (X : Scheme) (Z : PreClos X) (γ : Z.cov.J)
    [DecidableEq Z.indnumb]
    [Fintype Z.indnumb]
    [(i : Z.indnumb →₀ ℤ) → Decidable (i ∈ Set.range (ρNatToInt Z.indnumb))] : Scheme :=
  BlMu (A := Z.cov.obj γ) (ideal_loc X Z γ)

def open_pair (X : Scheme) (Z : PreClos X) (γ δ : Z.cov.J) : Scheme :=
  pullback (Z.cov.map γ) (Z.cov.map δ)

def Proj_loc_pair (X : Scheme) (Z : PreClos X) (γ δ : Z.cov.J)
  [DecidableEq Z.indnumb]
    [Fintype Z.indnumb]
  [(i : Z.indnumb →₀ ℤ) → Decidable (i ∈ Set.range (ρNatToInt Z.indnumb))] : Scheme :=
    pullback (pullback.fst (Z.cov.map γ) (Z.cov.map δ))

def Proj_loc_pair_mor (X : Scheme) (Z : PreClos X) (γ δ : Z.cov.J)
    [DecidableEq Z.indnumb]
    [Fintype Z.indnumb]
    [(i : Z.indnumb →₀ ℤ) → Decidable (i ∈ Set.range (ρNatToInt Z.indnumb))] :
  Proj_loc_pair X Z γ δ ⟶ open_pair X Z γ δ :=
    pullback.fst (pullback.fst (Z.cov.map γ) (Z.cov.map δ))
      (Proj_loc X Z γ ↘ Spec (Z.cov.obj γ))
\end{lstlisting}

Here \texttt{Proj\_loc} is $B_\gamma$, \texttt{open\_pair} is
$U_{\gamma\delta}$, \texttt{Proj\_loc\_pair} is $V_{\gamma\delta}$, and
\texttt{Proj\_loc\_pair\_mor} is $\pi_{\gamma\delta}$. The inclusion
$f_{\gamma\delta}$ is the other projection of the same fibre product.

\begin{proposition}\label{prop:swap}
There is exactly one isomorphism
\[
  t_{\gamma\delta} : V_{\gamma\delta} \xrightarrow{\ \sim\ } V_{\delta\gamma}
\]
over $U_{\gamma\delta}$.
\end{proposition}

This is the heart of the globalisation, and it is again an application of
affine terminality, now on the overlap. Both $V_{\gamma\delta}$ and
$V_{\delta\gamma}$ sit over $U_{\gamma\delta}$. Locally on $U_{\gamma\delta}$,
say on an affine open $\Spec C \subseteq U_{\gamma\delta}$, Proposition
\ref{prop:bc} identifies
\[
  V_{\gamma\delta} \times_{U_{\gamma\delta}} \Spec C
  \;\cong\; \BlM\bigl( I_{i\gamma} C \bigr),
  \qquad
  V_{\delta\gamma} \times_{U_{\gamma\delta}} \Spec C
  \;\cong\; \BlM\bigl( I_{i\delta} C \bigr).
\]
The two centres appearing on the right are equal, because both are the
restriction to $\Spec C$ of the global centre $Z$: the composites
$\Spec C \to U_\gamma \to X$ and $\Spec C \to U_\delta \to X$ agree, being the
two projections out of a fibre product over $X$. So the two sides are the
affine blowup of one and the same centre, and Theorem \ref{thm:aff-main} gives
a unique isomorphism between them over $\Spec C$.

Uniqueness at each $\Spec C$ is what allows these local isomorphisms to be
glued into $t_{\gamma\delta}$, and it also makes $t_{\gamma\delta}$ itself
unique. The statement proved in Lean records both the map and its restriction
behaviour.

\leansrc{Project/Blowups/BlowupExists.lean}{1963--1973}
\begin{lstlisting}
def Proj_loc_pair_lemm (X : Scheme) (Z : PreClos X) (γ δ : Z.cov.J)
  [DecidableEq Z.indnumb]
  [Fintype Z.indnumb]
  [(i : Z.indnumb →₀ ℤ) → Decidable (i ∈ Set.range (ρNatToInt Z.indnumb))] :
  ∃! (f : (Proj_loc_pair X Z γ δ) ⟶ (Proj_loc_pair X Z δ γ)),

  (∃ (pf : Scheme.Hom.IsOver f (open_pair X Z γ δ)),
    ∀ (C : CommRingCat) (i : Spec C ⟶ open_pair X Z γ δ) [IsOpenImmersion i],
      restrictToOpen f i =
      (pullbackSymmetry _ _).hom ≫ Proj_loc_pair_open_φ X Z γ δ C i ≫
      (pullbackSymmetry _ _).hom) := by
\end{lstlisting}

\noindent
The morphism $t_{\gamma\delta}$ itself is then extracted from that unique
existence statement.

\leansrc{Project/Blowups/BlowupExists.lean}{2259--2264}
\begin{lstlisting}
def Proj_loc_pair_swap (X : Scheme) (Z : PreClos X) (γ δ : Z.cov.J)
  [DecidableEq Z.indnumb]
  [Fintype Z.indnumb]
  [(i : Z.indnumb →₀ ℤ) → Decidable (i ∈ Set.range (ρNatToInt Z.indnumb))] :
    (Proj_loc_pair X Z γ δ) ⟶ (Proj_loc_pair X Z δ γ) :=
  Classical.choose (Proj_loc_pair_lemm X Z γ δ)
\end{lstlisting}

\subsection{Triple overlaps and the cocycle}\label{ssec:cocycle}

To glue the $B_\gamma$ along the $t_{\gamma\delta}$ we must supply the triple
overlap data and verify the cocycle condition. Let
\[
  P_{ijk} \;:=\; V_{ij} \times_{B_i} V_{ik},
\]
the preimage in $B_i$ of $U_i \cap U_j \cap U_k$. The required map is
\[
  t'_{ijk} : P_{ijk} \longrightarrow P_{jki},
\]
compatible with $t_{ij}$, and the cocycle condition is
\[
  t'_{ijk} \circ t'_{jki} \circ t'_{kij} = \mathrm{id}.
\]

The usual construction of $t'_{ijk}$ is topological: one checks that the image
of the triple overlap under $t_{ij}$ lies inside the double overlap
$V_{jk}$, and restricts. Here this is not needed, and the reason is worth
spelling out because it is the cleanest illustration of how the universal
properties organise the proof.

Recall that $V_{jk}$ was defined as a fibre product,
$V_{jk} = B_j \times_{\Spec A_j} U_{jk}$, and that $U_{jk}$ is itself a fibre
product $U_j \times_X U_k$. So to give a map $P_{ijk} \to V_{jk}$ it is enough
to give three maps
\[
  P_{ijk} \to B_j, \qquad P_{ijk} \to \Spec A_j, \qquad P_{ijk} \to \Spec A_k,
\]
agreeing over $\Spec A_j$ and over $X$ respectively. All three are available.
The first is
\[
  P_{ijk} \xrightarrow{\ \mathrm{pr}_1\ } V_{ij}
  \xrightarrow{\ t_{ij}\ } V_{ji} \hookrightarrow B_j,
\]
the second is its composite with $B_j \to \Spec A_j$, and the third is
\[
  P_{ijk} \xrightarrow{\ \mathrm{pr}_2\ } V_{ik}
  \xrightarrow{\ t_{ik}\ } V_{ki} \hookrightarrow B_k \to \Spec A_k .
\]
The compatibility over $X$ holds because $t_{ij}$ and $t_{ik}$ are morphisms
over the respective overlaps, so both composites $P_{ijk} \to X$ agree with the
structure morphism. So $t'_{ijk}$ is built by the universal property of fibre
products alone, with no reference to underlying topological spaces.

\leansrc{Project/Blowups/BlowupExists.lean}{2433--2446}
\begin{lstlisting}
noncomputable def Proj_loc_pair_t' (X : Scheme) (Z : PreClos X) (i j k : Z.cov.J)
    [DecidableEq Z.indnumb] [Fintype Z.indnumb]
    [(i : Z.indnumb →₀ ℤ) → Decidable (i ∈ Set.range (ρNatToInt Z.indnumb))] :
    pullback (Proj_loc_pair_incl X Z i j) (Proj_loc_pair_incl X Z i k) ⟶
      pullback (Proj_loc_pair_incl X Z j k) (Proj_loc_pair_incl X Z j i) :=
  pullback.lift
    (pullback.lift
      (pullback.lift
        (pullback.fst (Proj_loc_pair_incl X Z i j) (Proj_loc_pair_incl X Z i k) ≫
          Proj_loc_pair_swap X Z i j ≫ Proj_loc_pair_incl X Z j i ≫
          (Proj_loc X Z j ↘ Spec (Z.cov.obj j)))
        (pullback.snd (Proj_loc_pair_incl X Z i j) (Proj_loc_pair_incl X Z i k) ≫
          Proj_loc_pair_swap X Z i k ≫ Proj_loc_pair_incl X Z k i ≫
          (Proj_loc X Z k ↘ Spec (Z.cov.obj k)))
\end{lstlisting}

The three nested \texttt{pullback.lift}s are exactly the three maps above,
assembled into $V_{jk} = B_j \times_{\Spec A_j} (U_j \times_X U_k)$ from the
inside out.

For the cocycle, note that all the inclusions $f_{ij} : V_{ij} \to B_i$ are
open immersions, hence monomorphisms, so $P_{ijk} \to B_i$ is a monomorphism
too and it suffices to prove
\[
  \bigl( t'_{ijk} \circ t'_{jki} \circ t'_{kij} \bigr) \circ m_{ijk}
  \;=\; m_{ijk}, \qquad m_{ijk} : P_{ijk} \hookrightarrow B_i .
\]
Both sides are morphisms $P_{ijk} \to B_i = \BlM\bigl((I_{\bullet i})\bigr)$
from an open subscheme of $B_i$. By Section \ref{ssec:cars} the centre is
Cartier on that open, so Proposition \ref{prop:aff-uniq} applies: it is enough
to check the two morphisms agree after composing with
$B_i \to \Spec A_i$. Tracing through the definitions, that computation
telescopes and reduces to the single identity
\[
  t_{ij} \circ t_{ji} = \mathrm{id}_{V_{ij}} ,
\]
which is proved the same way, by the same uniqueness statement, using that a
double swap induces the identity on $U_{ij}$ since it is the composite of the
two symmetry isomorphisms of a fibre product.

So the entire triple overlap analysis reduces to one instance of affine
uniqueness. No intersections of opens are ever computed.

\subsection{Gluing and the global universal property}\label{ssec:global}

We can now assemble the gluing datum. Recall what has to be supplied: objects
$B_\gamma$, opens $V_{\gamma\delta} \hookrightarrow B_\gamma$ with
$V_{\gamma\gamma} = B_\gamma$, transition isomorphisms
$t_{\gamma\delta} : V_{\gamma\delta} \to V_{\delta\gamma}$ with
$t_{\gamma\gamma} = \mathrm{id}$, triple overlap maps
$t'_{\gamma\delta\epsilon}$ compatible with the $t_{\gamma\delta}$, and the
cocycle identity. We have produced all of these:
\begin{itemize}
\item $B_\gamma$ and $V_{\gamma\delta}$ in Section \ref{ssec:overlaps};
\item $t_{\gamma\delta}$ in Proposition \ref{prop:swap};
\item $t'_{\gamma\delta\epsilon}$ and the cocycle in Section
      \ref{ssec:cocycle}.
\end{itemize}
The two degenerate conditions are worth a comment. That $f_{\gamma\gamma}$ is
an isomorphism holds because $U_{\gamma\gamma} = U_\gamma \times_X U_\gamma$
equals $U_\gamma$, since $u_\gamma$ is a monomorphism. That
$t_{\gamma\gamma} = \mathrm{id}$ then follows from the uniqueness clause of
Proposition \ref{prop:swap}, since the identity is a morphism over
$U_{\gamma\gamma}$ and there is only one such. Both are places where an
informal account would not pause, and both need the monomorphism property of
open immersions to make the two available structure morphisms on
$V_{\gamma\gamma}$ agree.

The data $\bigl( B_\gamma, V_{\gamma\delta}, f_{\gamma\delta}, t_{\gamma\delta},
t'_{\gamma\delta\epsilon} \bigr)$ is therefore a gluing datum in the sense of
\cite{stacks-project}, and we define $\Bl_Z X$ to be the resulting scheme, with its
morphism to $X$ induced by the $B_\gamma \to \Spec A_\gamma \to X$. The
canonical maps $\iota_\gamma : B_\gamma \to \Bl_Z X$ are open immersions whose
images cover, and they satisfy
\[
  \iota_\gamma \circ f_{\gamma\delta}
  \;=\; \iota_\delta \circ f_{\delta\gamma} \circ t_{\gamma\delta}
\]
on $V_{\gamma\delta}$, which is the identity used in the gluing arguments
below. The glued scheme is the following, \texttt{PreBlGlob} being the gluing
datum whose fields are the constructions above.

\leansrc{Project/Blowups/BlowupExists.lean}{2586--2589}
\begin{lstlisting}
abbrev BlGlob (Z : PreClos X) [DecidableEq Z.indnumb]
  [Fintype Z.indnumb]
  [(i : Z.indnumb →₀ ℤ) → Decidable (i ∈ Set.range (ρNatToInt Z.indnumb))] : Scheme :=
  Scheme.GlueData.glued (PreBlGlob Z)
\end{lstlisting}

\begin{proposition}\label{prop:glob-cartier}
The centre pulled back to $\Bl_Z X$ is Cartier.
\end{proposition}

The presentation is assembled from the affine ones. Its cover is the family of
all dilatation charts of all the local blowups,
\[
  \bigl\{ \Spec A_\gamma[F] \;:\; \gamma \in J, \ F \in \Pres(I_{\bullet\gamma})
  \bigr\},
\]
indexed by the pairs $(\gamma, F)$, each mapping into $\Bl_Z X$ as the
composite
\[
  \Spec A_\gamma[F] \hookrightarrow B_\gamma
  \xrightarrow{\ \iota_\gamma\ } \Bl_Z X,
\]
a composite of open immersions and hence an open immersion. These charts cover,
since the $\iota_\gamma$ cover $\Bl_Z X$ and the dilatation charts cover each
$B_\gamma$.

The ideals are the ones supplied by Proposition \ref{prop:blm-cartier}, namely
$I_{i\gamma} \cdot A_\gamma[F]$, and they are principal on non-zero-divisors
for the same reason as there. What has to be checked in addition is that these
ideals really do present the pulled back centre, that is, that
\[
  \Spec\bigl( A_\gamma[F] / I_{i\gamma}A_\gamma[F] \bigr)
  \;\cong\;
  Z_i \times_X \Spec A_\gamma[F] ,
\]
compatibly with the structure morphisms. This follows by composing the
corresponding isomorphism for the presented centre $Z$ over $U_\gamma$ with the
base change of $Z_i \times_X U_\gamma$ along
$\Spec A_\gamma[F] \to \Spec A_\gamma$, using that a fibre product taken in two
stages agrees with the one taken in a single stage. In the formalisation this
composite is written out explicitly, and it is where most of the length of the
construction sits.

The Cartier conditions themselves are the affine ones, unchanged. This gives
Theorem \ref{thm:blglob-cartier}.

\leansrc{Project/Blowups/BlowupExists.lean}{3294--3300}
\begin{lstlisting}
def BlGlobPreClos (Z : PreClos X)
    [DecidableEq Z.indnumb] [Fintype Z.indnumb]
    [(i : Z.indnumb →₀ ℤ) → Decidable (i ∈ Set.range (ρNatToInt Z.indnumb))] :
    PreClos (BlGlob Z) where
  indnumb := Z.indnumb
  subscheme i := pullback (Z.subscheme i ↘ X) ((BlGlob Z) ↘ X)
  over i := ⟨pullback.snd _ _⟩
\end{lstlisting}

The choice of cover here is the whole content. A cover of $\Bl_Z X$ by
arbitrary affine opens would not do, since the pulled back centre is invertible
but not necessarily principal on such an open. It is because the Cartier
condition of Definition \ref{def:cartier} quantifies existentially over
presentations that we are free to use the cover we built the scheme from, on
which principality holds by construction.

For the universal property, let $T \to X$ be such that the centre becomes
Cartier on $T$. Set $T_\gamma := T \times_X U_\gamma$. On $T_\gamma$ the centre
$I_{\bullet\gamma}$ becomes Cartier, by restricting the given presentation, so
Theorem \ref{thm:aff-main} gives a unique
\[
  \varphi_\gamma : T_\gamma \longrightarrow B_\gamma
\]
over $\Spec A_\gamma$. Composing with $B_\gamma \to \Bl_Z X$ gives maps
$T_\gamma \to \Bl_Z X$ over $X$, and these glue: on
$T_\gamma \times_T T_\delta$ both restrictions are morphisms over the overlap
into a blowup of a Cartier centre, so they agree by Proposition
\ref{prop:aff-uniq}. Note again that the compatibility is not computed but
deduced from uniqueness. This produces $\varphi : T \to \Bl_Z X$ over $X$, and
gives Theorem \ref{thm:main-pre}'s existence half.

Uniqueness is the same argument with the maps given rather than constructed. If
$\varphi, \varphi'$ are two morphisms $T \to \Bl_Z X$ over $X$, they agree
after restriction to each $T_\gamma$ by Proposition \ref{prop:aff-uniq}, hence
agree.

\subsection{Independence of the presentation}\label{ssec:rel}

It remains to see that $\Bl_Z X$ depends only on the centre and not on the
presentation, which is what allows Theorem \ref{thm:main} to be stated for
$\cZ \in \Clos(X)$.

\begin{proposition}\label{prop:rel}
If $Z'$ and $Z''$ are equivalent presented centres, then there is exactly one
isomorphism $\Bl_{Z'} X \cong \Bl_{Z''} X$ over $X$.
\end{proposition}

\leansrc{Project/Blowups/BlowupExists.lean}{3354--3361}
\begin{lstlisting}
lemma PreProjBlowup_rel (Z' Z'' : PreClos X) (eq : Quotient.mk' Z' = Quotient.mk' Z'')
    [DecidableEq Z'.indnumb]
    [Fintype Z'.indnumb]
    [Fintype Z''.indnumb]
    [DecidableEq Z''.indnumb]
    [(i : Z'.indnumb →₀ ℤ) → Decidable (i ∈ Set.range ⇑(ρNatToInt Z'.indnumb))]
    [(i : Z''.indnumb →₀ ℤ) → Decidable (i ∈ Set.range ⇑(ρNatToInt Z''.indnumb))] :
  ∃! (e : BlGlob Z' ≅ BlGlob Z''), Scheme.Hom.IsOver e.hom X := by
\end{lstlisting}

Once Theorem \ref{thm:blglob-cartier} is available this is formal, and it is
the same two-sided argument as Proposition \ref{prop:bc}. Since $Z'$ and $Z''$
represent the same element of $\Clos(X)$, the centre that becomes Cartier on
$\Bl_{Z'}X$ is the same element of $\Clos$ as the one that becomes
Cartier on $\Bl_{Z''}X$. So Theorem \ref{thm:main-pre} applied to $Z''$ with
$T = \Bl_{Z'}X$ gives a unique map one way, the symmetric application gives a
unique map the other way, and both composites are identities by uniqueness.

Here Theorem \ref{thm:main-pre} is used with its hypothesis stated on the class:
$f^*[Z]$ is Cartier (\texttt{PreProjBlowup\_UnivProp\_\allowbreak existence\_Cars},
\texttt{PreProjBlowup\_UnivProp\_\allowbreak unicity\_Cars}). Since $[Z'] = [Z'']$,
this hypothesis is the same for $Z'$ and $Z''$.

We may therefore define, for a centre $\cZ \in \Clos(X)$,
\[
  \Bl_\cZ X \;:=\; \Bl_{Z} X \quad \text{for any presentation } Z \text{ of } \cZ,
\]
well defined up to unique isomorphism over $X$, and Theorem \ref{thm:main}
follows from Theorem \ref{thm:main-pre} by rewriting the Cartier hypothesis
along the equality $[\,Z\,] = \cZ$. This completes the proof.

\subsection{Dependencies}\label{ssec:deps}

The logical structure of the proof is the following, read left to right.

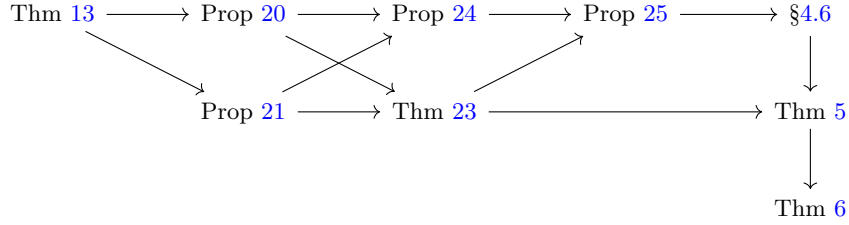
\[
\begin{tikzcd}[column sep=1.1cm, row sep=0.75cm]
  \text{\small Thm \ref{thm:dil-univ}}
    \arrow[r] \arrow[dr]
  & \text{\small Prop \ref{prop:blm-cartier}}
    \arrow[r] \arrow[dr]
  & \text{\small Prop \ref{prop:bc}}
    \arrow[r]
  & \text{\small Prop \ref{prop:swap}}
    \arrow[r]
  & \text{\small \S\ref{ssec:cocycle}}
    \arrow[d] \\
  & \text{\small Prop \ref{prop:aff-uniq}}
    \arrow[r] \arrow[ur]
  & \text{\small Thm \ref{thm:aff-main}}
    \arrow[ur] \arrow[rr]
  & & \text{\small Thm \ref{thm:main-pre}}
    \arrow[d] \\
  & & & & \text{\small Thm \ref{thm:main}}
    \arrow[d] \\
  & & & \text{\small Prop \ref{prop:cartier-mathlib}}
    \arrow[r]
  & \text{\small Cor \ref{cor:mathlib}}
\end{tikzcd}
\]

Everything originates in the universal property of the dilatation. The affine
terminality statement, Theorem \ref{thm:aff-main}, is used four separate times
downstream: to build the base change isomorphism, to build the overlap
comparison, to close the cocycle, and to glue the local maps into the global
one. Proposition \ref{prop:aff-uniq}, the uniqueness half, is used more often
still, since it is what makes each of those constructions canonical and hence
compatible with the others.

Note that the diagram has no cycles and no step that depends on the global
object being already constructed. This is not automatic. A common way to
organise a gluing argument is to construct the object first and verify its
properties afterwards, which works informally because one can speak of the
object before knowing that it is well defined. Here every arrow points forward,
and the global scheme $\Bl_Z X$ is named only after the gluing datum has been
completely verified.

\section{Correspondence}\label{ssec:correspondence}

For reference, the statements of this paper and the declarations proving them.

\begin{center}
\renewcommand{\arraystretch}{1.25}
\begin{tabular}{@{}lll@{}}
\hline
Statement & Lean declaration & File \\
\hline
Def.\ \ref{def:preclos} & \texttt{PreClos} & \texttt{PreClosAndClos.lean} \\
Def.\ \ref{def:rel} & \texttt{rel}, \texttt{Clos} & \texttt{PreClosAndClos.lean} \\
Def.\ \ref{def:cartier} & \texttt{IsPreCars}, \texttt{IsCars} & \texttt{PreClosAndClos.lean} \\
\S\ref{ssec:pullback} & \texttt{pullback\_PreClos}, \texttt{pullback\_Clos} & \texttt{PreClosAndClos.lean} \\
Def.\ \ref{def:multicenter} & \texttt{Multicenter}, \texttt{LargeIdeal} & \texttt{Multicenter.lean} \\
Def.\ \ref{def:dilatation} & \texttt{PreDil}, \texttt{r}, \texttt{Dilatation} & \texttt{Multicenter.lean} \\
Lem.\ \ref{lem:dil-nzd} & \texttt{nonzerodiv\_image} & \texttt{Multicenter.lean} \\
Lem.\ \ref{lem:dil-gen} & \texttt{image\_elem\_LargeIdeal\_equal} & \texttt{Multicenter.lean} \\
Thm.\ \ref{thm:dil-univ} & \texttt{desc}, \texttt{lemma\_exists\_unique\_morphism} & \texttt{Multicenter.lean} \\
Cor.\ \ref{cor:dil-trivial} & \texttt{ofFamilyIso} & \texttt{Multicenter.lean} \\
Prop.\ \ref{prop:dil-funct} & \texttt{functo\_dila\_alg} & \texttt{Multicenter.lean} \\
Def.\ \ref{def:pres} & \texttt{Mu} & \texttt{Bl.lean} \\
Prop.\ \ref{prop:dil-is-potion} & \texttt{clo\_mu\_mor}, \texttt{Mu\_mor\_iso} & \texttt{Bl.lean} \\
Def.\ \ref{def:blm} & \texttt{BlMu} & \texttt{Bl.lean} \\
Prop.\ \ref{prop:blm-cartier} & \texttt{BlMuPreClos\_IsPreCars}, \texttt{blowups\_Cars} & \texttt{Cars.lean} \\
Prop.\ \ref{prop:aff-uniq} & \texttt{ProjBlowup\_UnivProp\_unicity\_affine} & \texttt{UniqueBlowup.lean} \\
Prop.\ \ref{prop:aff-exist} & \texttt{ProjBlowup\_UnivProp\_existence\_affine} & \texttt{BlowupExists.lean} \\
Thm.\ \ref{thm:aff-main} & \texttt{ProjBlowup\_is\_conceptual\_blowups\_affine} & \texttt{BlowupExists.lean} \\
Prop.\ \ref{prop:bc} & \texttt{base\_change\_Bl\_open} & \texttt{BlowupExists.lean} \\
\S\ref{ssec:overlaps} & \texttt{Proj\_loc}, \texttt{Proj\_loc\_pair} & \texttt{BlowupExists.lean} \\
Prop.\ \ref{prop:swap} & \texttt{Proj\_loc\_pair\_lemm} & \texttt{BlowupExists.lean} \\
\S\ref{ssec:cocycle} & \texttt{Proj\_loc\_pair\_t'}, \texttt{PreBlGlob} & \texttt{BlowupExists.lean} \\
\S\ref{ssec:global} & \texttt{BlGlob} & \texttt{BlowupExists.lean} \\
Prop.\ \ref{prop:glob-cartier} & \texttt{BlGlobPreClos}, \texttt{BlGlob\_IsCars} & \texttt{BlowupExists.lean} \\
Thm.\ \ref{thm:main-pre} & \texttt{PreProjBlowup\_UnivProp} & \texttt{BlowupExists.lean} \\
Prop.\ \ref{prop:rel} & \texttt{PreProjBlowup\_rel} & \texttt{BlowupExists.lean} \\
Thm.\ \ref{thm:main} & \texttt{GlobalBlowup\_UnivProp} & \texttt{BlowupExists.lean} \\
Prop.\ \ref{prop:cartier-mathlib} & \texttt{isCars\_pullback\_toClos\_iff} & \texttt{IdealSheaf.lean} \\
Cor.\ \ref{cor:mathlib} & \texttt{exists\_blowup\_isEffectiveCartier} & \texttt{IdealSheaf.lean} \\
\hline
\end{tabular}
\end{center}

\noindent
Files in \texttt{Project/Blowups/} except \texttt{Multicenter.lean}, which is in
\texttt{Project/Dilatation/}.

\bmhead{Acknowledgements}
We are very grateful to the referee for a careful reading of the manuscript and
for many valuable comments. We especially thank the referee for pointing out
some inconsistencies in the mathematical presentation of an earlier version, and for suggestions
that improved the formalization and its exposition.

\end{document}